\documentclass[12pt]{article}
\usepackage{MyStyle-2}

\usepackage{graphicx} 
\usepackage{amsfonts}
\usepackage{titlesec}
\usepackage{booktabs}
\usepackage{tabularx}
\usepackage{placeins} 
\usepackage[
  backend=biber,
  style=alphabetic,
  backref=true
]{biblatex}
\titleformat{\subsubsection}[runin]
  {\normalfont\normalsize\bfseries}
  {\thesubsubsection}
  {1em}
  {}
\usepackage{amssymb,amsmath,amsthm}
\usepackage{xcolor}
\usepackage{quiver}

\theoremstyle{plain}

\title{Genus one mirror symmetry and surfaces with trivial first Chern class}
\author{Gabriel Abánades-Joglar, Dennis Eriksson }
\date{\today}

\author{
Gabriel Abánades-Joglar
\and
Dennis Eriksson 
}
\begin{document}
\maketitle

\begin{abstract}
The main goal of this paper is to study the extension of genus-one BCOV type mirror symmetry for smooth projective varieties with vanishing $c_1(X)$. We show that the naive generalization of the expectation that one can encode genus one Gromov-Witten invariants in the BCOV invariant of the mirror family cannot hold in this setup. We show this by completing the study of these invariants for surfaces. Progress had been made recently by Oberdieck's proof of the Klemm-Mariño conjecture for Enriques surfaces, which can be combined with computations of equivariant analytic torsion by Yoshikawa. We also perform the analogous study for bielliptic surfaces. 

Both of these cases present a sign issue not present in the known trivial canonical bundle cases, which we show implies that the mirror symmetry predictions cannot hold for varieties with torsion canonical bundle, even up to sign.

\end{abstract}

\setcounter{tocdepth}{2}\tableofcontents 

\section{Introduction}

\subsection{Main results} In the early $90$'s the study of string Theory led physicists to formulate surprising predictions regarding the enumerative geometry of Calabi--Yau manifolds. These motivated, in particular, the development and construction of Gromov--Witten invariants. In genus one, these predictions are stated by Bershadsky--Cecotti--Ooguri--Vafa in physical terms for Calabi--Yau threefolds in \cite{BCOV_K-S_TheoryGravity} and \cite{bershadsky1993holomorphic}. First instances of these predictions were both interpreted mathematically and proven by Fang--Lu--Yoshikawa for the quintic threefold in \cite{Fang-Lu-Yoshikawatorsionthreefolds} and by Eriksson--Freixas i Montplet--Mourougane for Calabi--Yau hypersurfaces in projective space in \cite{eriksson2022genus}. All of these deal with smooth projective varieties with trivial canonical bundle. 

In the setting of varieties with torsion canonical bundle, Klemm--Mariño study string theory in the local Enriques threefold, given by the total space of the canonical bundle of an Enriques surface. Using the language of BCOV they make interesting predictions regarding the Gromov--Witten invariants of an Enriques surface. This conjecture was recently proven by Oberdieck in \cite{Oberdieck2025KMformula}. At the same time, there are computations of what should be the $B$-model correlation function in the context of K3 surfaces with involution \cite{YoshikawaK3inventiones}. Both computations involve a certain orthogonal modular form given by the Borcherds $\Phi$-function. It is tempting to see this as a mirror symmetry statement, and part of the purpose of this article is to interpret these things in the context of \cite{eriksson2022genus}. 

One of the main takeaways of our work can be interpreted as the following ''meta-theorem''
\begin{maintheorem}\label{intro: main theorem}
Let $X$ be a projective manifold with $c_1(X)=0$. 
    \begin{itemize}
        \item The Deligne--Riemann--Roch formalism of \cite{eriksson2022genus} extends to the setting of torsion canonical bundle and one can define a function $F_1^B$.
        
        \item For surfaces one can find a mirror family $X^\vee\to (\D^\ast)^r$ with a change of variables $q: \D^r \to \D^r$, called the mirror map, such that the relationship 
        \begin{displaymath} F_1^A(\tau) = -F_1^B(q(s)).
        \end{displaymath}  
        \item The naive generalization of the conjecture, namely that,
        \begin{displaymath}
            F_1^A=\pm F_1^B
        \end{displaymath} cannot hold.
    \end{itemize}
\end{maintheorem}
We dedicate the following paragraphs to render precise and contextualize the contents of this theorem.

\subsection{$B$-side potential for trivial Chern classes}

In this section we work with projective complex manifolds $X$ with $c_1(X)=0$ of dimension $n$. Equivalently, they can be seen as complex manifolds admitting a Ricci-flat metric or as complex manifolds with torsion canonical bundle. The setting at hand is  intended to extend previous discussions having actually trivial canonical bundle. 

\subsubsection{} As mentioned above we are interested in studying connections between complex geometric data and enumerative data. The enumerative invariants are encoded in the so--called A--model potential, which is a function on the complexified Kähler cone $\mathcal{K}_\C(X)=H_\R^{1,1}/H^{1,1}_\Z+iK(X)$ where $K(X)$ is the usual Kähler cone.
The potential is defined as follows:
\begin{equation}\label{Intro:def of F1A general}
F_1^A(\tau)=\frac{-1}{24}\left(\int_Xc_{n-1}(X)\wedge 2\pi i \tau\right) +\sum_{\beta\in H_2(X,\Z)}N_{1,\beta}q^\beta
\end{equation}
where $\tau$ is a complexified Kähler class and $q^\beta=e^{2\pi i \langle\beta,\tau\rangle}$. Here $N_{1,\beta}$ is the unpointed genus one count of curves in a given curve class $\beta$, formally defined in terms of Gromov--Witten invariants (see \cite{BeherendGWinvariants}). These counts give numbers and don't need insertions exactly because of the assumption that $c_1(X)=0$, see for example \cite[\textsection 1.1]{eriksson2022genus}.

\subsubsection{} The second object we are interested in is the so-called BCOV invariant. It is a real-valued invariant of a Calabi--Yau manifold,  introduced first by Fang--Lu--Yoshikawa in \cite{Fanglu} for threefolds and in the general setting of trivial canonical bundle by the second author together with Freixas i Montplet and Mourougane 
\cite{eriksson2021bcov}. Finally, it was extended to pluricanonical sections by Zhang \cite{Yeping.zhang_extensionBCOV}. We refer the reader to section \ref{Sect: BCOV invariant c1=0} for a review and precise discussion in the torsion canonical bundle case. Given $X\to S$ a family of compact Kähler manifolds of dimension $n$ with vanishing canonical class this invariant is defined as 
$$\tau_{BCOV}(X/S)=\frac{A(X/S)}{B(X/S)}\prod_{p,q}\det({\Delta^{p,q}_{\overline{\partial}}})^{(-1)^{p+q}pq}, $$
where $\det(\Delta^{p,q}_{\overline{\partial}})$ denotes the zeta regularized determinant of the $\overline{\partial}$ operator acting on $(p,q)$ forms and $A,B$ are certain renormalizing factors to render the left-hand-side independent of auxiliary metrics. This specific product of determinants was introduced in the physics literature by Bershadsky--Cecotti--Ooguri--Vafa (BCOV) in \cite{BCOV_K-S_TheoryGravity}. 

\subsubsection{} Suppose now that $S = (\D^\ast)^r$, the punctured multidisc, with coordinates $s_i$. We suppose that monodromies of $X \to S$ around each puncture are unipotent. Extending previous discussions, our results in \ref{subsection: Riemann--roch formulation of BCOV} can be stated as the following theorem: 
\begin{theorem}\label{intro: thm about extending DRR formalism} Suppose the function $\log \tau_{BCOV}$ is $L^2_{loc}$. Then there exists a holomorphic factorization
\begin{equation}\label{eq: holo factorization intro}
    \tau_{BCOV}(X/S)^{6} =C\cdot |\exp((-1)^n 24F_1^B)| g(s)
\end{equation}
Such that the following properties hold:
\begin{itemize}
    \item $F_1^B$ is uniquely determined by the expansion 
    \begin{displaymath}
        F_1^B = \sum_{i=1}^r \alpha_i  \log s_i+\text{ holo in $s$ with no constant term}
    \end{displaymath}
    where $s=(s_1,\ldots, s_r)$ is the coordinate of the (punctured) multidisc and all $\alpha_i \in \R$.
    \item $g(s)$ is an explicit combination of $L^2$-style norms of sections of Hodge-theoretic bundles determined by the limit mixed Hodge structure.
    \item  $C$ is an undetermined constant. 
\end{itemize}
\end{theorem}
Notice that the above theorem does not require us to be in the large complex structure limit in order for $F_1^B$ to be defined. The only necessary geometric condition is that the monodromy is unipotent.
The $L^2_{loc}$ assumption is verified in the algebraic setting. It could moreover be lifted if the Deligne--Riemann--Roch isomorphism (\cite{eriksson2025deligneriemannrochisomorphism}) is proven to be an isometry, which is expected. 

\subsubsection{} In the setting of mirror symmetry, given a projective manifold $X$ with $K_X=\Ocal_X$ one expects, see for example \cite[Section 6, Conjecture 2]{eriksson2022genus} and the surrounding discussion, that there exists a large complex structure degeneration $X^{\vee} \to (\D^\ast)^r$ together with a mirror coordinate change $s \mapsto q(s)$ that satisfies

\begin{equation}\label{eq: MS for hypersurfaces intro}
    F_{1,X}^A(\tau)=F_{1,X^\vee}^B(q(s)).
\end{equation}

This actually has been shown to hold in different cases in \cite{Fang-Lu-Yoshikawatorsionthreefolds,eriksson2022genus,eriksson2025genusmirrorsymmetryintersection}.

Thanks to theorem \ref{intro: thm about extending DRR formalism}  it makes sense to wonder if in the torsion canonical bundle setting there is any relation between $F_1^A$ and $F_1^B$.

\subsection{The case of Calabi--Yau surfaces}
We will henceforth suppose that $X\to S$ is a family of  surfaces with torsion canonical bundle. Turns out, by the Kodaira classification of compact algebraic surfaces, that any such surface is either an abelian surface, a $K3$-surface, an Enriques surface or a bielliptic surface.

\subsubsection{} Let us begin by saying a few words about the function $g$ appearing in Theorem \ref{intro: thm about extending DRR formalism} when we restrict to surfaces. As we mentioned before this function is given as a combination of $L^2$ norms of sections of two line bundles which are studied in section \ref{Sect: BCOV invariant c1=0}. The first line bundle, given in general in equation \eqref{eq: def BCOV bundle}, reduces to $\lambda(\Ocal)^{-2}$, where $\lambda$ denotes the determinant of cohomology thanks to the isometric properties of Serre duality 
(\cite[Prop. 3.4]{Fang-Lu-Yoshikawatorsionthreefolds}). The second bundle we study is $\left(f_\ast (K_{X/S}^{\otimes m})\right) ^{\frac{\chi(X_s)}{12m}}$ where $m$ is the torsion of the canonical bundle of a fiber in the family.

 Below we give a list with the description of both bundles in the four cases at hand.

\begin{table}[h]\label{table: DRR surfaces}
\centering
\renewcommand{\arraystretch}{1.25}
\begin{tabular}{|c|c|c|c|c|}
\hline
Type & $m$ & $\chi_{\mathrm{top}}$ &
Fiberwise Dolbeault determinant line 
& $(f_\ast K_{X/S}^{\otimes m})^{\frac{\chi}{12m}}$
\\
\hline
K3 & 1 & 24 &
$\det H^{2,0}(X)^{2}$
&
$(f_\ast K_{X/S})^{2}$
\\
\hline
Abelian surface & 1 & 0 &
$\Ocal_S$
&
$\Ocal_S$
\\
\hline
Enriques & 2 & 12 &
$\Ocal_S$
&
$(f_\ast K_{X/S}^{2})^{1/2}$
\\
\hline
Bielliptic (hyperelliptic) & 2,3,4,6 & 0 &
$\det H^{0,1}(X)^{2}$
&
$\Ocal_S$
\\
\hline
\end{tabular}

\end{table}
\FloatBarrier

\subsubsection{} The mirror symmetry setting also requires one to consider a large complex structure limit which correspond to very degenerate cusps of the moduli space. Normally these degenerations have been identified with maximally unipotent degenerations (MUM) in the case of threefolds (\cite{Morrisonguide}). In the case of Calabi--Yau hypersurfaces in $\P^n$, inspired by Deligne  in \cite{eriksson2022genus}, the first author and collaborators showed that Hodge--Tate degenerations work as a choice of large complex structure. 

In the case  of torsion canonical bundle there is an associated family $Y/S$ of varieties with trivial canonical bundle, together with a fiber-wise finite degree cyclic covering $Y_s\to X_s$, given by a group $G$. The formalism of distinguished sections require the family $Y/S$ to be unipotent in middle cohomology.  In order to include the case of bielliptic surfaces into this formalism, no stronger conditions can be expected to hold. 

\subsubsection{The cusp problem}\label{subsubsect:cuspproblem intro} Already in the case of moduli spaces of Enriques surfaces, there are two natural candidates for such degenerations in this setup, namely the level 1  and the level 2 cusps of moduli of Enriques surfaces. Only the level 1 cusp admits an identifiable mirror symmetry phenomenon, meaning that simply a degeneration alone cannot determine a mirror family. The case of bielliptic surfaces present even more possible choices of cusps, corresponding to the cusps of various modular curves. Our study of analytic torsion shows that only one of these cusps allows for an interpretation of mirror symmetry phenomena.

\subsubsection{} With the following theorem we present our suggestion on how to choose a cusp in the torsion canonical bundle case: 
\begin{theorem}\label{thm: MS for surfaces thm intro}
    Let $X$ be any Calabi--Yau surface. Then: 
    \begin{enumerate}
        \item There exists a mirror family $\mathcal{X}^\vee \to (\D^\ast)^r$, such that:
        \begin{enumerate} 
        \item The variation of Hodge structures associated with  $H^k(\mathcal{X}_t)$ has a limit mixed Hodge structure which is Hodge--Tate.
        
        \item At the level of coarse moduli the cusp represented by the limit of the covering family $\mathcal{Y}^\vee_t$ corresponds to a (possibly) singular variety with an action of $G$ without fixed points.
        
    \end{enumerate}
        \item The exists a mirror map $q:(\D^\ast)^r\to (\D^
    \ast)^r$ given in terms of canonical periods.
    \end{enumerate}
    such that 
    $$F_1^A( \tau)=-F_1^B(q(s)).$$
Here  $\tau\in \mathcal{K}_\C(X)$ is an element in the complexified Kähler cone. 
\end{theorem}
Below we discuss the various cases that need to be treated individually. 
\subsubsection{The case of K3 and abelian surfaces} In this setting the theorem follows from the fact that the usual (non-reduced) Gromov-Witten theory is trivial, i.e. $GW_1(\beta)=0$ for all $\beta \in H_2$ and it follows that we have $F_1^A=0$. This follows from \cite{kiem2013localizing} and it is well known in general. On the B-side this is \cite[Prop 5.12]{eriksson2021bcov}. In particular the distinguished sections are the unit of $\mathcal{O}_X$ for the abelian case, and an integral holomorphic $2$ form for the $K3$. In these cases $F_1^B=0$. Hence the second part of the theorem is simply the equality $0=-0$.

\subsubsection{The case of Enriques.} 
Our contribution to this case focuses on constructing a family of Enriques surfaces using a toroidal compactification of a specific period domain. As we mentioned there is a choice of cusp to be made between level 1 or level 2. Let us now be more explicit about the mirror symmetry phenomena around the level 1 cusp.

In \cite{Oberdieck2025KMformula} Oberdieck obtains that
$$F_1^A=\log \left(\Phi ^{-\frac{1}{8}}\right),$$
where the Borcherds $\Phi$-function is expanded around the level 1 cusp. On the other hand, Yoshikawa's result in \cite{YoshikawaK3inventiones} implies that 
\begin{displaymath} \tau_{BCOV}=C\cdot ||\Phi||_{\text{Pet}}^{\frac{1}{2}}.
\end{displaymath}
Our theorem in this case reduces to proving all of these expressions have an interpretation in terms of distinguished sections and their canonical periods.

In this setting we remark that the chosen cusp have a property as mentioned Theorem \ref{thm: MS for surfaces thm intro} by the extension results of \cite{alexeev2023compactificationsmodulispacesk3}. Their results are for special $1$-dimensional families degenerating into the cusps in the coarse moduli space.  We have not given a general moduli interpretation of this along other limits going towards this cusp. For more details around this discussion see subsection \ref{subsect: Enriques dist sections}.

\subsubsection{The case of bielliptic surfaces.} The main body of this work consists of the computation of the $B$-side function for bielliptic surfaces. These are surfaces given as a quotient of a product of two elliptic curves $E\times F$ by a specific group action $G=\Z/N\Z\times \Z/K\Z$. We recall their geometry in Section \ref{sect: geom biell}. Their complex moduli are essentially that of the modular data of level $N$, i.e. the curve $X_1(N)$.

The key computation that we carry out is that of the $BCOV$ invariant, in Theorem \ref{thm:computation of BCOV invariant}. We show that:
\begin{equation}\label{eq: bcov bielliptic intro}
\tau_{BCOV}=\left(2\im(\tau')\right)^2\prod_{j|N}\left(j^2|\Delta(j\tau')|^{\frac{1}{6}}\right)^{\lambda(j)}
\end{equation}
where $\tau '$ is a specific transformation of the period of $F$ adapted to a specific cusp in the modular curve $X_1(N)$ that ensures the third property in \ref{thm: MS for surfaces thm intro} and $\Delta$ is the discriminant modular form.

One can explicitly check from our computations that if we don't apply the mirror map (equivalently, if we choose another cusp) $F_1^B$ and $F_1^A$ are not clearly related, in particular they don't have the same leading term. 

In the last section we carry out the computation of genus--one Gromov--Witten invariants of bielliptic surfaces. This is done by using the degeneration formula in the style of the computations by Blomme in \cite{blomme2025gromov} where he computes the generating functions $F_g^A$ for $g\geq 2$. We obtain exactly the following formula:
\begin{equation}\label{intro: GW of biell}
    F_1^A=-\left(\sum_{j|N} \lambda(j)\log(\eta(q^j))\right)
\end{equation}
where the  $\lambda(j)$ are the same as the exponents appearing in the BCOV invariant and $\eta=\Delta^{\frac{1}{24}}$ is the Dedekind eta.

An interesting observation appearing in this case is that the Hodge theoretical part of the formalism (including the distinguished sections) is governed by $H^1 (X)$ rather than the middle cohomology $H^2$, as is usually expected from manifolds with trivial canonical bundle since they satisfy local Torelli.

\subsection{Towards mirror symmetry for trivial Chern class.}\label{subsect: general conjecture }

As we recalled after equation \eqref{eq: MS for hypersurfaces intro} the results of \cite{Fang-Lu-Yoshikawatorsionthreefolds,eriksson2022genus} for Calabi--Yau hypersurfaces say that, after applying the mirror map, one has $F_1^A=F_1^B$.

\subsubsection{Positive results.} 
 We can produce new Calabi--Yau manifolds by taking products of them. It is natural to investigate the compatibility of the mirror symmetry conjecture in this context.  
 \begin{theorem}
     Suppose $X,Y$ are Calabi--Yau manifolds of dimensions $n,m$ with trivial canonical bundle. If the mirror symmetry expectation $F_1^A=F_1^B$ is true for both $X$ and $Y$, then it also is for $X \times Y$.
 \end{theorem}
The mirror family we consider is the product of mirror families. That is, if $\mathcal{X}^\vee/(\D^\ast)^{r_1}$ and $\mathcal{Y}^{\vee}/(\D^\ast)^{r_2}$ are the mirror families we compute $F_1^B$ for the product family $\mathcal{X}^{\vee}\times \mathcal{Y}^{\vee}/(\D^\ast)^{r_1+r_2}$.

The key computation is the following two identifications in Theorem \ref{thm: M-S products}
\begin{equation}\label{eq:F1Bintro}
    F_1^B(\mathcal{X}_s ^{\vee}\times \mathcal{Y}_t^{\vee })=(-1)^m \chi(\mathcal{Y}_t^{\vee }) F_1^B(\mathcal{X}_s^{\vee })+ (-1)^n \chi(\mathcal{X}_t^{\vee }) F_1^B(\mathcal{Y}_s^{\vee }).
\end{equation}
and
\begin{equation}\label{eq:F1Aintro}
   F_1^A(X\times Y)=\chi (Y)F_1^A(X)+\chi(X)F_1^A (Y).
\end{equation}
The first one amounts to a computation of the BCOV invariant, cf. Lemma \ref{lemma:Factor T for product} of a product and the compatibility of products and sections, cf. Prop. \ref{prop:iso lines and sections of product}. The second one amounts to a computation of the relevant Gromov--Witten invariants, coming from the work of Behrend \cite{BeherendProduct}. 

The theorem then follows from the additional information that for mirror pairs we have the expectation $\chi(X)=(-1)^n \chi(\mathcal{X}_s^\vee).$

Notice that here the result is in the trivial canonical bundle case, this allows for new examples of traditional enumerative mirror symmetry in genus one. 

\subsubsection{Negative results.}\label{subsubsect:negativeresults intro}As a contrast to what we showed above, we present now an example that shows that the naive generalization $$F_1^A= \pm F_1^B$$ cannot hold if we consider the whole class of varieties with torsion canonical bundle. This means that the there isn't a general conjectural image in this setup, even though all the invariants generalize with good properties and there are semi-positive results for surfaces.

\begin{theorem}
    Let $S$ be an Enriques surface and $Q$ be a quintic threefold. Then if $X = Q \times S$ we do not have $F_1^A=\pm F_1^B$.
\end{theorem}
The counter-example will also rely on \eqref{eq:F1Bintro} and \eqref{eq:F1Aintro} which are valid also for torsion canonical bundles. By the main result of \cite{Fang-Lu-Yoshikawatorsionthreefolds} (or \cite{eriksson2022genus} for the arbitrary hypersurface) we know that $Q$ admits a mirror family $Q^\vee\to \D^\ast$, by our theorem \ref{thm:mirror theorem Enriques} we know that $S$ admits a mirror family $S^\vee\to (\D^\ast)^{10}$. It is natural then, to consider $X^\vee=Q^\vee \times S^\vee\to (\D^\ast)^{11}$ as the mirror family of $X=Q\times S$. 

Since $\chi(S)=\chi(S_s^\vee)=12$ and $\chi(Q)=-\chi(Q_t^\vee)=-200$ we have that $$F_1^B(X^\vee)=-200 F_1^B(S^\vee_s)+  12 F_1^B(Q^\vee_t)$$ 
and
$$F_1^A(X)= -200 F_1^A(S)+12 F_1^A(Q).$$
While looking tantalizingly similar, since $F_1^A(S)=-F_1^B(S^\vee_s)$ and $F_1^A(Q)=F_1^B(Q_t^\vee)$, these do not match up to sign. 

\subsubsection{} It is tempting to speculate on a possibly correct formulation of the mirror symmetry statement in the setting of projective manifolds with trivial Chern class. We plan to investigate this phenomenon in future work. 

\subsubsection*{AI disclosure}
The mathematical results present are all the work of the authors. We used LLMs to verify spelling and some presentation aspects of the text. Part of the classical results for the K3 period domains presented in subsection \ref{subsect: moduli polarized K3} we learned in part by some added discussions with a model.

While generating a revision of typos the model mentioned a statement in the lines of our proposition \ref{prop: GW are zero for b_1 nonzero} but the present proof is our own. In the preliminary version checked by the LLM we had a different approach to lemma \ref{lema:logDRR is psh} by assuming only $L^1_{loc}$ and using results of Xia. However, the LLM pointed out that this reference was wrong in this generality, so we took the present approach through a more classical method by assuming $L^2_{loc}$.

\subsubsection*{Acknowledgments} To be added.

\newpage

\section{The BCOV invariant for trivial Chern class manifolds}\label{Sect: BCOV invariant c1=0}
This section sets up the general framework needed to pose the refined BCOV conjectures of Eriksson--Freixas i Montplet--Mourougane for projective manifolds such that $c_1=0$ in $H^2(X,\Q)$, equivalently $K_X^{\otimes m}=0$ for some $m\in \mathbb{N}$. The first thing we do is present the BCOV invariant in this context. This is a real valued $C^\infty$ function given as the quotient of two metrics on a line bundle. These metrics appear as renormalizations of the Quillen and $L^2$-metrics so that they become independent of the Kähler structure. 

\subsubsection{}The construction of the  renormalized $L^2$-metric is that of \cite{eriksson2021bcov}, since it works for arbitrary compact Kähler manifolds. We present the main properties of this metric in Prop.\ref{prop:renormalized L^2 norm}. In particular, in \ref{item3:renormalized L^2 norm} we correct a small mistake of \cite[Remark 5.5]{eriksson2021bcov} concerning the expression of the factor $B(X)$ when using Poincar\'e duality. This is due to the fact that we choose to use the Arakelov-theoretic convention to choose our Kähler forms. This is to make it so Poincar\'e duality is a morphism of Hodge structures $$H^k(X,\Z)\otimes H^{2n-k}(X,\Z)\to H^{2n}(X,\Z)=\Z(-n).$$ 

The renormalization of the Quillen metric is done in the style of  Zhang's construction \cite{Yeping.zhang_extensionBCOV} where the BCOV is extended to a more general setting. We present here a direct approach in the style of \cite{Fang-Lu-Yoshikawatorsionthreefolds} to also explicitly show where certain $2\pi$-factors appear because we consider Arakelov-theoretic conventions for the Kähler forms. The main properties of the renormalized Quillen metric are stated in Prop.\ref{prop:QuillenBCOV}.

\subsubsection{} After this is done we study the Hodge-theoretical properties of the BCOV bundle. This is the parallel discussion to \cite[Subsect. 6.2]{eriksson2022genus} but done for torsion canonical bundle now. The main idea here is to introduce what we mean by strongly unipotent degenerations in this setting. We show that for these degenerations $f:X\to S$ both $f_\ast K_{X/S}^m$ and $\lambda_{BCOV}$ admit distinguished sections $\eta_0,\eta_{BCOV}$ that trivialise them. 

In contrast to what happens for usual Calabi-Yau manifolds, where both trivializations come from Hodge theory only, in the present setting the construction of $\eta_0$ requires certain geometric information coming from the degree $m$ cover of $X$ associated to the triviality of $K_X^{\otimes m}$. This formalism of distinguished sections allows to consider certain holomorphic periods which we interpret as (certain pieces of) a mirror map since in the classical case these correspond to canonical coordinates in the tangent space of the moduli space of Complex structures of $X$.

Once these sections are constructed we can compare (certain powers of) their $L^2$ and Quillen-BCOV metrics as in \cite{eriksson2022genus}. By this procedure we obtain a function $|DRR|$ which conjecturally in the Calabi--Yau case realizes mirror symmetry (see \cite[Section 6, Conjecture 2]{eriksson2022genus}), as we mention in the introduction there is no clear way to relate $|DRR|$ and $F_1^A$ in the general settin of torsion canonical bundle (see \ref{subsubsect:negativeresults intro}).

\medskip

\subsection{BCOV torsion and invariant for trivial Chern class}
In this section we consider $X$ a compact K\"ahler manifold with trivial (in $H^2(X,\Q)$) Chern class. By the Beauville--Bogomolov decomposition it follows that the canonical bundle is in fact torsion. In this and greater generality, a BCOV invariant was introduced by Zhang in \cite{Yeping.zhang_extensionBCOV}, generalizing constructions from \cite{eriksson2021bcov} and \cite{Fang-Lu-Yoshikawatorsionthreefolds}. In this section, to clarify some normalizations and for the convenience of the reader, we introduce the BCOV invariant in the case the canonical bundle being torsion. The specific construction is also important for the geometric reformulation in terms of Grothendieck--Riemann--Roch.

\subsubsection{} Recall first that the BCOV torsion is defined as \begin{equation}\label{eq:BCOV-torsion}\mathcal{T}_{BCOV}=\prod_{p,q} \left(\det\left(\Delta^{p,q}_{\overline{\partial}}\right)\right)^{(-1)^{p+q}pq},
 \end{equation} 
where the determinant of the Laplacian is the zeta-regularized determinant. This is a priori dependent on the implicit K\"ahler metric. In order to renormalize it to the effect that is independent, we follow the techniques from \cite{eriksson2021bcov} and  \cite{Fang-Lu-Yoshikawatorsionthreefolds}. 
To this end, we recall the definition of the BCOV line given in the previous cited works as: 
\begin{equation}\label{eq: def BCOV bundle}
    \lambda_{BCOV}=\bigotimes_{p=0}^n\lambda(\Omega^p_{X/S})^{(-1)^pp}
\end{equation}

The BCOV torsion naturally appears as the quotient of the Quillen metric $h_Q$ and $h_{L^2}$ on $\lambda_{BCOV}$. We will renomalize these metrics, to obtain two metrics $h_{Q,BCOV}$ and $h_{L^2,BCOV}$,  that are independent of the Kähler metric. We will then define the BCOV invariant to be \begin{equation}\label{eq:defbcov}
\tau_{BCOV}=\frac{h_{Q,BCOV}}{h_{L^2,BCOV}}=\frac{A(X)}{B(X)}\mathcal{T}_{BCOV}
\end{equation}
where $A,B$ are the normalizing factors we have to find and $\mathcal{T}_{BCOV}$ is the BCOV torsion. 

\subsubsection{The renormalized $L^2$-norm} We recall for the convenience of the reader, and future reference, the renormalized $L^2$-norm in \cite{eriksson2021bcov}. We begin with the Arakelov-theoretic $L^2$ norm. That is, if $h(-,-)$ is the Kähler metric on $X$ we take the Kähler form locally as $\omega=\frac{i}{2\pi}\sum h(\frac{\partial}{\partial z_i},\frac{\partial}{\partial\overline{z}_j})dz_i\wedge d\overline{z}_j$. The $L^2$ norm is then 
$$h_{L^2}(\alpha,\beta)=\int_Xh(\alpha,\beta)\frac{\omega^n}{n!}$$ In order to make it Kähler independent, we define 
\begin{equation}\label{def:B-term}B(X/S)=\prod_{k=0}^{2n}\vol_{L_2}(H^k(X_s,\Z),\omega)^{(-1)^{k+1}\frac{k}{2}}
\end{equation}
Where $\vol(H^k(X,\Z),\omega)$ denotes the square of the covolume of the lattice determined by the image of the cohomology $H^k(X,\Z)$ inside $H^k(X,\R)$ computed with respect to the $L^2$ inner product associated to $\omega$.
\begin{definition}
    The $BCOV$ normalized $L^2$ norm is the norm on $\lambda_{BCOV}$ given by:
    $$h_{L^{2},BCOV}=B(X)h_{L^2}$$
\end{definition}

For the convenience of the reader we recall the main properties of the $L^2-BCOV$ norm in the following proposition which can be extracted, with the appropriate care about the factors of $2\pi$ that we explain in the proof,  from \cite[Proposition 4.2 ,Remark 5.5 (ii),Proposition 5.6]{eriksson2021bcov}:
\begin{proposition}\label{prop:renormalized L^2 norm}
    Suppose that $X \to S$ is a proper submersion of relative dimension $n$ of complex manifolds. The following properties hold:
    \begin{enumerate}
        \item\label{item1:renormalized L^2 norm} The $BCOV$-normalized $L^2$ norm is independent of the Kähler forms used for defining the $L^2$-norms.
        \item\label{item2:renormalized L^2 norm}If $X$ admits a global Kähler form then the factor $B(X/S)$ is a locally constant function on $S$.
        \item \label{item3:renormalized L^2 norm} One has $$\vol(H^k(X_s,\Z),\omega)\vol(H^{2n-k}(X_s,\Z),\omega)=\frac{1}{(2\pi)^{2b_kn}}$$ and so the following equality is satisfied 
        $$B(X/S)=\prod_{k=0}^{n-1}\vol(H^{k}(X_s,\Z),\omega)^{(-1)^k(n-k)}\cdot (2\pi)^{-n\left(\sum_{k=0}^{n-1}(-1)^kkb_k+(-1)^nb_n\frac{n}{2}\right)}\cdot (2\pi)^{\chi(X)n^2}.$$
    \end{enumerate}
\end{proposition}
\begin{proof}
    Items \ref{item1:renormalized L^2 norm} and \ref{item2:renormalized L^2 norm} are proven in \cite{eriksson2021bcov}. There, a similar statement to \ref{item3:renormalized L^2 norm} is given but there is a discrepancy of factors of $2\pi$. We introduce the following notation $v_k:=\vol(H^k(X_s,\Z),\omega)$

    The equality of expressions for $B(X/S)$ follows from the equality $v_k\cdot v_{2n-k}=\frac{1}{(2\pi)^{2b_kn}}$. Usually one would expect $v_k\cdot v_{2n-k}=1$ but we are using a renormalized Kähler form and thus our Serre and Poincaré duality pairings are also renormalized.
    In particular we have (see \cite{Fang-Lu-Yoshikawatorsionthreefolds}):
    \begin{equation}\label{Poincar\'e duality}
    \begin{aligned}
        H^k(X,\Z)\times H^{2n-k}(X,\Z)  \to &Z(-n) \\
        (\alpha,\beta)  \xmapsto{\hspace{1cm}} \left(\frac{1}{i2\pi}\right)^n&\int_X\alpha\wedge \beta
    \end{aligned}
    \end{equation}
    
    This map is induced by the direct sums of the Serre duality pairings $$SD(\omega_1,\omega_2)=\frac{1}{(i2\pi)^n}\int_X\omega_1\wedge \omega_2$$ for $\omega_1\in H^{p,q}$, $\omega_2\in H^{n-p,n-q}$.
    Serre duality induces (see for example \cite{Fang-Lu-Yoshikawatorsionthreefolds}) the isometry $\det(H^{p,q})\otimes \det(H^{n-p,n-q})=(\C,|\cdot |)$. On the other hand the Poincar\'e duality pairing \ref{Poincar\'e duality} gives $\det(H^k(X,\Z))\otimes \det(H^{2n-k}(X,\Z)\otimes \Z(n))\cong \Z$. Which implies $$\det(H^k(X,\Z))\otimes \det(H^{2n-k}(X,\Z))\cong \Z(-b_{2n-k}n).$$ 
    Since this isomorphism is induced by Serre duality it is an isometry and thus gives the desired equality for volumes.
\end{proof}

\subsubsection{The renormalized Quillen metric}

We will modify the construction of \cite{Fang-Lu-Yoshikawatorsionthreefolds} to renormalize the Quillen metric on the BCOV line bundle.

For this, we start by observing that the line bundle $f_\ast(K_{X/S}^{\otimes m})$ has a metric only dependent on the complex structure, given by 
\begin{equation}\label{def:L2m-norm}||\eta||^2_{L^{\frac{2}{m}}}(s):=\left( \frac{1}{(2\pi)^n}\int_{X_s} |\eta_s\wedge \overline{\eta}_s |^{\frac{1}{m}}\right)^m.
\end{equation}
Where $\eta_s$ is a section of the trivial bundle $K_{X_s}^{\otimes m}$.

Moreover, since \eqref{eq:evalmapcanonical} is an isomorphism, we can always, after possibly shrinking $S$, choose a global nowhere vanishing section $\eta\in H^0(S,f_\ast K^m_{X/S})$. We likewise think of it as a section of $K_{X/S}^{m}$, necessarily non-vanishing on the fibers.

We refer to 
$$\omega_{WP}(X/S)=-\frac{1}{m}dd^c\log \|\eta\|^2_{L^{\frac2m}} $$
as the Weil--Petersson form.

With these definitions in mind, we introduce the following analogue of the factor $A$ of \cite[Definition 4.1]{Fang-Lu-Yoshikawatorsionthreefolds}, which in fact is independent of $\eta$:
    \begin{equation}\label{def:A factor}
    A(X/S)=\vol(X_s,\omega_s)^{\frac{\chi(X_s)}{12}}\exp\left(\frac{-1}{12}\int_{X/S}\log\left(\frac{|\eta\wedge \overline{\eta}|^{\frac{1}{m}}}{\left(\frac{((2\pi)\omega)^n}{n!}\right)}\frac{\vol(X_s,\omega_s)}{\left(||\eta||^2_{L^{\frac{2}{m}}}\right)^{\frac{1}{m}}}\right)c_n(T_X,h)\right)
    \end{equation}
This is a function of $s \in S$. Notice also that the volume factors outside and inside the integral cancel out but, as we will see shortly in Proposition \ref{prop:QuillenBCOV}, it is very convenient to write $A$ in the above way.

\begin{definition}
    The BCOV-normalized Quillen norm on $\lambda_{BCOV}$ is defined as 
\begin{equation}
    h_{Q,BCOV}:=A(X/S)\cdot h_Q,
\end{equation}
where $h_Q=\mathcal{T}_{BCOV} \cdot h_{L^2}$ is the usual Quillen metric.
\end{definition}

We are ready to state the following proposition:
\begin{proposition}\label{prop:QuillenBCOV}
    The following properties hold:
    \begin{enumerate}
        \item \label{item1:QuillenBCOV} The first Chern form of the Quillen-BCOV metric is equal to  $\frac{\chi(X_s)}{12m}c_1(f_\ast K^{\otimes m}_{X/S},h_{L^{\frac2m}})=\frac{\chi(X_s)}{12}\omega_{WP}(X/S)$.
        \item \label{item2:QuillenBCOV} The BCOV-normalized Quillen norm is independent of the auxiliary K\"ahler metric. 
        \item \label{item3:QuillenBCOV} If the K\"ahler metric is Ricci flat, the term inside the logarithm in \eqref{def:A factor} is 1, and hence 
        $A(X/S) = \vol(X_s, \omega)^{\frac{\chi(X_s)}{12}}.$
    \end{enumerate}
\end{proposition}

\begin{proof}
Property \ref{item1:QuillenBCOV} actually implies \ref{item2:QuillenBCOV}, since the curvature only depends on the complex structure. This follows by interpolating two different K\"ahler metrics over a constant family parametrized by $\mathbb{P}^1$. The details are as in \cite[Theorem 4.10]{Fang-Lu-Yoshikawatorsionthreefolds}.

To establish \ref{item1:QuillenBCOV} , we will use the curvature formula of Bismut-Gillet-Soul\'e \cite{BGS1988}, in the form of \cite[Theorem 3.5]{Fang-Lu-Yoshikawatorsionthreefolds}. For our family $f:X\to S$ the formula reads 

\begin{equation}\label{eq:curvatureBCOV}
c_1(\lambda_{BCOV},h_Q)=\frac{-1}{12}\int_{X/S}c_1(T_{X/S},h)c_n(T_{X/S},h).
\end{equation}
    We will compute the curvature of the Quillen metric on $\lambda_{BCOV}$ without the factor $A(X/S)$.

    We know that we can compute the first Chern form globally, in terms of a global non-vanishing section $\eta$ of $K_{X/S}^m$:
    \begin{equation}\label{eq:ChernformTX}-c_1(T_{X/S}, h) = c_1(K_{X/S},h)=\frac{1}{m}c_1(K_{X/S}^m,h)=\frac{-1}{m}dd^c \log(h(\eta,\eta))=-dd^c \log\left(\frac{|\eta\wedge \overline{\eta}|^{\frac{1}{m}}}{\left(\frac{\omega^n}{n!}\right)}\right).
    \end{equation}
    Inserting \eqref{eq:ChernformTX} into \eqref{eq:curvatureBCOV}, and using that Chern forms are $d$ and $d^c$-closed, immediately leads to:

\begin{equation}\label{eq: sum of chern forms}
    \begin{split}
    \int_{X/S}c_1(T_{X/S},
    h)c_n(T_{X/S},h) &=   \int_{X/S}dd^c\frac{1}{m}\log(||\eta||^2_{L^{\frac{2}{m}}})c_n(T_X,h)-dd^c \log(A(X/S)^{12}) \\
    & =  \frac{1}{m}dd^c \int_{X/S} \log(||\eta||^2_{L^{\frac{2}{m}}}) c_n(T_X, h)-dd^c \log(A(X/S)^{12})\\
    & = -\frac{\chi(X_s)}{m} c_1(f_\ast K_{X/S}^m, h_{L^{2/m}})- dd^c \log(A(X/S)^{12}).
    \end{split}   
\end{equation}

In the last equality we used the Gauss--Bonnet formula. Dividing by $-12$ proves the first point.

To obtain part \ref{item3:QuillenBCOV} we can assume that we are in the case where $S$ is a point. Here  Ricci flatness simply means $c_1(K_X,h_{K_X})=0$. By \eqref{eq:ChernformTX} this means that $\log\left(\frac{|\eta\wedge \overline{\eta}|^{\frac{1}{m}}}{\frac{\omega^n}{n!}}\right)$ is a pluriharmonic function on a compact manifold, and thus constant. We compute the constant value $\mu$ inside the $\log$  by integrating. We find that $\mu \int_X\frac{\omega^n}{n!}=\int_X|\eta\wedge\overline{\eta}|^{\frac{1}{m}}$ which gives $\mu=\frac{{ (2\pi)^n(||\eta||^2_{L^{2/m}}})^{\frac1m}}{vol(X,\omega)}$ from this part \ref{item3:QuillenBCOV} follows.

\end{proof}

\subsection{Riemann--Roch formulation - geometric and analytic}\label{subsection: Riemann--roch formulation of BCOV}

For the below, we refer to \cite[\textsection 6]{eriksson2022genus} for a more detailed account. That discussion involves the case of trivial canonical bundle, and we modify the needed pieces here accordingly.  Let $f: X \to S$ be a family of projective manifolds, with connected fibers. We suppose there is a non-zero integer $m$ such that the canonical bundle is $m$-torsion, i.e. there is an isomorphism $K_{X_s}^{m}\simeq \Ocal_{X_s}$.  In this case, the evaluation morphism  
\begin{equation}\label{eq:evalmapcanonical}
    f^\ast f_\ast (K_{X/S}^{m}) \to K_{X/S}^{m}
\end{equation}
is necessarily an isomorphism. Also, for any $p$, introduce the line bundles 
\begin{equation}\label{eq:detcohomology}
    \lambda(\Omega^p_{X/S}) = \det Rf_*(\Omega_{X/S}^p).
\end{equation}
Fiberwise, they are isomorphic to the alternating product of determinants of Dolbeault cohomology groups $\bigotimes_q \det H^{q}(X, \Omega_{X_s}^p)^{(-1)^q}$. 

In this setting, the functorial Riemann--Roch theorem \cite{eriksson2025deligneriemannrochisomorphism} implies that : 
\begin{equation}\label{eq: DRR for BCOV bundle}
    \bigotimes \lambda(\Omega_{X/S}^p)^{(-1)^p p} \simeq (f_\ast K_{X/S}^m)^{\frac{\chi}{12m}}
\end{equation}
where $\chi$ is the topological Euler characteristic of the fibers of $f$. It is an isomorphism of $\Q$-line bundles. 

We are interested in studying the local behavior of the family around certain points of the base. From the point of view of Mirror Symmetry we want the degenerations to be maximal. In order to define what we mean my maximally degenerate around a point, we assume that $S=(\mathbb{D}^\ast)^ d$.
\subsubsection{}\label{subsubsect: strongly unipotent}
Recall also that we say that the family is strongly unipotent close to a point 0, if $S \simeq (\D^\ast)^d$ and 
\begin{itemize}
    \item for any $i$, the local monodromies of the local systems $R^i f_\ast \C$ are all unipotent.
    \item for any $i$, the limit mixed Hodge structures $H^i_{lim}$ are Hodge--Tate, meaning they only have $(p,p)$ parts.
\end{itemize}
In this setting we review and adapt to our context the conjectures in \cite{eriksson2022genus} relating mirror symmetry and distinguished sections of the BCOV line.

\subsubsection{}\label{subsubsuect dist section}
Firstly, let us briefly recall what happens in the classical Calabi-Yau case. In this context the degeneration being strongly unipotent gives that $\mathcal{F}^p/\mathcal{F}^{p+1}\cong \mathcal{W}_{2p}/\mathcal{W}_{2p-2}$ as vector bundles. This in turn, since the monodromy sends $W_k$ to $W_{k-2}$, gives that both bundles are trivial. Using the lattice structure coming from $H^k(X,\Z)$ and this trivialization, one can construct canonical sections unique up to sign (see the discussion in \cite[Section 6]{eriksson2022genus}, in particular Remark 6.8). 

During this text we will choose as canonical sections the ones in \textit{loc.cit.} multiplied by a power of $2\pi i$.  Due to the fact that we want the comparison morphisms between homology and cohomology to be morphisms of Hodge structures. In particular, the definition of distinguished section in the present context is the following (see \cite[Lemma 6.9]{eriksson2022genus}): 

\begin{definition}\label{def:Hodge Tate distinguished sections}
Let $(\mathcal{F}^\bullet,H^k)$ a geometric variation of Hodge structures with a Hodge-Tate limit for the monodromy weight filtration $W_\bullet$.

Let $\eta$ be a section of $\mathcal{F}^p\subset H^k$. We say that its class $\eta_{p,k-p}$ in $\mathcal{F}^p/\mathcal{F}^{p+1}$ is a distinguished section if 
    \begin{itemize}
        \item $\int_\gamma \eta=0$ for all $\gamma\in W_{-2p-2}H_k$
        \item $\frac{1}{(2\pi)^n}\int_\gamma \eta$ is constant for all $\gamma\in W_{-2p}H_k$
    \end{itemize}
    
\end{definition}
\subsubsection{}\label{distinguished section for 1 dimensional cohomology}
    For $h^{p,q}=1$ with $k<2n$ we have that the integral (with respect to the Twisted lattices) distinguished section is given by the condition 
\begin{equation}
    \frac{1}{(2\pi)^n}\int_{\gamma_0}\eta=1 \text{ for } \gamma_0 \text{ generator of } W_{-2n}H_k(Z,\Q)\cap H_k(X,\Z).
\end{equation}
Observe that for $2n$ the power of $2\pi$ is already implicit in our Kähler form, that is, the (Tate twisted) dual of the point class is $\frac{\omega^n}{n!}=:1_{n,n}$.

In the Calabi-Yau case, since $F^nH^n=H^0(X,K_X)$, we have a construction for canonical sections of the right hand side of DRR in \ref{eq: DRR for BCOV bundle} also. We describe our approach to the case of torsion canonical bundle now. If the degeneration of the family $X/S$ is still Hodge-Tate then \cite[Lemma 6.4, Remark 6.5]{eriksson2022genus} still apply and we have a well defined notion of distinguished sections of the BCOV line. However, it is not clear what a distinguished section of $f_\ast K^{m}_{X/S}$ should be. 

\subsubsection{}The construction explained above gives sections that are canonical up to $\operatorname{GL}(n,\Z)$ by using the fact that the graded objects of the weight filtration carry a lattice. We now explain how to construct sections of certain pieces of variations of Hodge structures. These sections will only be canonical up to an element of $\operatorname{GL}(n,\C)$, our main goal is to construct them for $\mathcal{F}^n$,  the deepest piece of the Hodge filtration of the 
middle cohomology of a Calabi--Yau family that is not necessarily Hodge--Tate. We assume that the monodromy is unipotent.

The sections we construct will be the ones associated to the opposite filtration construction of  Fern\'andez--Pearlstein in \cite{OpositeFilFernandezPearlstein} (we call $S_k$ what they call $\Phi_k$). Given the limit mixed Hodge structure $(V,F_\infty,W(N))$ and its associated Deligne splitting 
$$V=\oplus_{r,s}I^{r,s},$$
one has 
$$F^p=\oplus_{r\geq p}I^{r,s},\;\; W_k=\oplus_{r+s\leq k}I^{r,s}.$$
Fern\'andez--Pearlstein consider the following increasing filtration, a priori different from the weight filtration:
$$S_k=\oplus_{r\leq k}I^{r,s}.$$

In the case where the MHS is Hodge--Tate, by definition, we only have $I^{r,r}$ appearing and one recovers $W_{2k}=S_k$. Among the many properties of the filtration $S$ studied in \cite{OpositeFilFernandezPearlstein} their Theorem 3.28 gives that this filtration extends, over a product of punctured disks, to a complex flat sub-vector bundle $(\Scal, \nabla) \subset (\mathcal{V}, \nabla^{GM})$ opposed to the Hodge filtration of the starting VHS. Namely, for any $p$, we have that the natural map 
\begin{displaymath}
    \Fcal^p \oplus \mathcal{S}_{p-1}\to V
\end{displaymath}
is an isomorphism. 

Moreover, for each element of the monodromy cone, we actually have: 
$$NS_k\subset S_{k-1}.$$
By a standard argument, as in the end of Lemma \ref{lemma:holoperiod} below, this gives that each monodromy logarithm around each puncture sends $\Scal_k$ to $\Scal_{k-1}$.
    If the monodromy is unipotent, this last condition means that the local system $\Gr_k^\Scal$ is trivial. As in the case of Hodge--Tate (compare with \cite[Eq. 6.5]{eriksson2022genus}), the fact that $\Scal$ is opposed to $\Fcal$ gives a canonical identification 
    \begin{equation}\label{eq:distinguishediso} \Gr^p_\Fcal=\Gr^\Scal_p.
    \end{equation}
    By the above remarks the right hand side is a trivial local system. 
    \begin{definition}
        We call any sections of $\Gr^p_{\Fcal}$, for a variation of Hodge structures as above distinguished, if they correspond via \eqref{eq:distinguishediso} to a flat, hence trivial, section of $\Gr_p^\Scal$.
    \end{definition}
    We remark moreover that while the limit mixed Hodge structure might depend on the choice of a parameter, the graded pieces $\Gr_p^\Scal$ are canonical for reasons analogous to that the graded pieces of the weight filtration are \cite[Paragraph after Eq. 6.15]{Schmid1973}.

\subsubsection{} 
For this setting, suppose that $f: X \to (\D^\ast)^r$ is a projective family of smooth manifolds of dimension $n$. Informally, the following lemma now states that distinguished sections are preserved under the Serre duality isomorphism $\det\left(\Gr^p_{\Fcal}H^k\right)\simeq \det\left(\Gr^{n-p}_{\Fcal}H^{2n-k}\right)^\vee.$ The statement can only be up to constant because already the definition of distinguished sections of $\det\Gr^p_{\mathcal{F}}$ is up to constants. 
\begin{lemma}
    Suppose $\eta$ is a distinguished section of $\det\left(\Gr^p_{\Fcal}H^k\right)$ and $\eta'$ is a distinguished section of $\det\left(\Gr^{n-p}_{\Fcal}H^{2n-k}\right)$ then the product of the norms $||\eta||_{L^2}||\eta'||_{L^2}$ is constant on $\left(\D^\ast\right)^r$, 
\end{lemma}
\begin{proof}
The result will follow from existence and the commutativity, up to constant, of the following diagram of flat vector bundles. 
\[\begin{tikzcd}
	{\mathcal{S}_p/\mathcal{S}_{p-1}} & {\mathcal{S}_p\cap\mathcal{F}^p} & {\mathcal{F}^p/\mathcal{F}^{p+1}} \\
	{\left(\mathcal{S}_{n-p}/\mathcal{S}_{n-p-1}\right)^\vee} & {\left(\mathcal{S}_{n-p}\cap\mathcal{F}^{n-p}\right)^{\vee}} & {\left(\mathcal{F}^{n-p}/ \mathcal{F}^{n-p+1}\right)^\vee}
	\arrow[from=1-2, to=1-1]
	\arrow[from=1-1, to=2-1]
	\arrow[from=1-2, to=1-3]
	\arrow[from=1-2, to=2-2]
	\arrow[from=1-3, to=2-3]
	\arrow[from=2-1, to=2-2]
	\arrow[from=2-3, to=2-2]
\end{tikzcd}\]
The horizontal arrows are isomorphism because $\mathcal{S}$ splits $\mathcal{F}$. The rightmost vertical arrow is the Serre duality isometry, the middle vertical arrow is Poincar\'e duality and the leftmost vertical arrow is the induced pairing on the graded pieces of the filtration $\mathcal{S}$. 

To justify first the commutativity, we notice that since Serre duality and Poincar\'e duality coincide up to a constant the right square commutes. The pairing on the quotient for the filtration $\mathcal{S}$ is defined by taking Poincar\'e duality on representatives. It passes to the quotient $\Gr_p^\Scal$ because the first Hodge--Riemann relation is preserved on the limit. Since the Poincar\'e pairing is a flat, also on the flat structures on the $\Gr_p^\Scal$, it sends flat  sections to flat sections, hence a distinguished basis to a distinguished basis of the two $\Gr_\Fcal$. The lemma is now an immediate consequence.
\end{proof}
\subsubsection{} The above arguments work for general families with unipotent monodromy, and will be utilized for families of bielliptic surfaces. The following lemma, and the discussion after it, shows that under certain conditions compatible with the Mirror Symmetry philosophy, one can actually find a section of $\mathcal{F}^n$ canonical up to sign.

\begin{lemma}\label{lemma: covering associated to torsion KX}
    Suppose that $K_{X/S}^m \simeq \Ocal_X$. Consider the associated order $m$-cyclic covering $Y \to X$, and denote by $h: Y\to X \to S$. Then: 
    \begin{itemize}
        \item $Y$ has trivial relative canonical bundle.
        \item We have an isomorphism
    $$ h_\ast K_{Y/S}^{m}\simeq f_\ast K _{X/S}^m.$$ 

    This isomorphism has norm $m^m$ when we take the $L^2$ and $L^{\frac{2}{m}}$ norms, respectively on each side. So we treat it as an isometry.
    \item If the family $Y/S$ has maximally unipotent monodromy at middle cohomology then $$h^{n,n}_{\lim}(Y)=1.$$
    \end{itemize}

\end{lemma}
\begin{proof}
    
The fact that $K_{X_s}^{\otimes m}$ has a nowhere vanishing global section $\eta$ gives that there exists a degree $m$ \'etale cover $\rho: Y_s\to X_s$ such that $Y_s$ has trivial canonical bundle. We also denote by $\rho: Y\to X$ the non-fiberwise covering. A fiberwise description of $Y\to X$ is the cyclic degree $m$ \'etale cover
$$Y_s=\{(x,\xi)\in X_s\times K_X|\xi\in (K_{X_s})_x \text{ and }\xi^{m}=1\}.$$  Since $\rho$ is \'etale we know $\rho^\ast (K_{X_s})=K_{Y_s}$. Thus, the section $\xi$ is a canonical section trivializing $K_{Y/S}$.
In algebraic terms (see for example \cite{VanishingEsnaultViehweg}) the cover is described as:
\begin{displaymath}
Y= \Spec_X \left(\bigoplus_{i=0}^{m-1} K_{X/S}^{-i}\right)
\end{displaymath}
This is an alternative way  of seeing that $\rho^{\ast}K_{X/S} = K_{Y/S}$ and finding the desired trivialization. The isomorphism between the push-forwards is given as follows:

We use that, since the fibers $Y_s$ are CY,  $h_\ast (K_{Y/S}^{\otimes m})=\left(h_\ast K_{Y/S}\right)^{\otimes m}$. Given $\alpha$ a local section of $h_\ast K_{Y/S}$ we have that $\alpha^{\otimes m}= g\rho^\ast \eta$ with $g$ a function on the base. We have that $\alpha^{\otimes m} \mapsto g\cdot \eta=\eta'$ gives the desired isomorphism. It is immediate, since $\rho$ is degree $m$ \'etale, that the $L^2$-norm of $\alpha ^{\otimes m}$ is $m^m$ times the $L^{\frac2m}$-norm of $\eta '$.

For the last part, we simply notice that if we have maximally unipotent monodromy, $N^{n}$ is non-trivial and induces an isomorphism $\Gr^p_F \Gr_{2n}^W H^{n}_{\lim}(Y_s)\simeq \Gr^{p-n}_F \Gr_0^W H^n_{\lim}(Y_s) $. Since $\Gr_0^W$ is pure of weight $0$, this means that $p=n.$ On the other hand, $\Gr_F^n H^n_{\lim}(Y_s)=F^nH^n_{\lim}(Y_s)$ is one-dimensional, because the general fiber is Calabi--Yau.

\end{proof}

\subsubsection{} Thus, from the perspective of Hodge theory, it is sensible to assume that our setting is the following: 
\begin{itemize}
    \item The family $f:X\to S$ has strongly  unipotent monodromy.
    \item The family $Y\to S$ associated to the covering of $X$ has unipotent monodromy in middle cohomology.
\end{itemize}
We assume this below. We have that $\mathcal{F}^n(Y)=\Gr_n^\mathcal{S}$ is one dimensional, so it admits a distinguished generator $\gamma_0^\vee$ unique up to constant by the flatness of the splitting.  We interpret this as the  dual of a certain cycle in homology. By construction, this cycle is killed by the  action on homology of any element $N$ of the monodromy cone  as $(\Gr_{n}^S)^\vee=S_{-n}H_{n}$ and we have $N(S_kH_n)\subset S_{k-1}H_n$ together with $S_{-n-1}H_n=0$. Here $S_{-k}$ denotes the dual of the filtration $S$.

\begin{definition}
    The distinguished section $\eta_0$ of $f_\ast K^{m}_{X/S}$ is the section  corresponding via $h_\ast K_{Y/S}=f_\ast K_{X/S}^{m}$ to $(\gamma_0)^\vee=\frac{\eta}{\frac{1}{(2\pi i)^n}\int_{\gamma_0}\eta}$ where $\eta$ is a nowhere vanishing section of $h_\ast K_{Y/S}$.
\end{definition}

We record the following lemma describing properties of the integral appearing in the above definition.

\begin{lemma}\label{lemma:holoperiod}
    Suppose that $S=(\D^\times)^d$ and that there is an extension of the family $X \to (\D^\times)^d$ to a projective $\widetilde{f}:\widetilde{X} \to \mathbb{D}^d$ with smooth total space. If $\eta$ extends as a holomorphic section of $\widetilde{f}_\ast K_{\widetilde{X}/\D^d}$, the period $s\mapsto\int_{\gamma_0(s)}\eta_s$ is a holomorphic function on the poly-disk $\mathbb{D}^d$. 
\end{lemma}
\begin{proof}
    A priori the map $s\mapsto \int_{\gamma_0(s)}\eta_s$ is multivalued. It is standard (see \cite[Equation 4.8]{Schmid1973}) that when one has a multivalued period map, the modification $s\mapsto e^{-
    \sum_{i=1}^d N_i \log s_i}\int_{\gamma_0}\eta_s$ becomes univalued (it is monodromy-invariant) and extends to all of $\mathbb{D}^d$. In this case, if we prove $N_i\gamma_0=0$ then this function coincides with $s\mapsto \int_{\gamma_0(s)}\eta_s$ so the period is univalued. 

    The fact that $N_i\gamma_0=0$ follows from the fact that the monodromy weight filtration $W(N)$ is independent of the element of the cone and hence so is $S$ (with fixed Hodge filtration of course). Since any $N$ is $(-1,-1)$ we have that for $N\in \{\sum_{i=1}^n\lambda_iN_i|\lambda_i>0\}$,  $NS_{l}\subset S_{l-1}$. Say $n=2$, then $(\frac{1}{2}N_1+N_2)S_l\subset S_{l-1}$ and $(N_1+N_2)S_{l}\subset S_{l-1}$. This implies $N_1S_{l}=2[(N_1+N_2)-(\frac{1}{2}N_1+N_2)]S_l\subset S_{l-1}$. The same argument gives $N_2S_{l}\subset S_{l-1}$ and by induction we have the desired result.  
    
    Now,  $\gamma_0^\vee$ trivializes the Deligne extension of $\mathcal{F}^n$ and $\eta$ is a section  $\widetilde{f}_\ast K_{\widetilde{X}/\D^d}$. By \cite[Theorem 2]{Semipostheorem}, and Remark 4.7 for the application to the analytic category, the natural identification outside the coordinate axes extends to all of $\D^d$. Hence $\eta=G(s)\gamma_0^\vee$ with $G(s)$ holomorphic. By construction $G(s)=\frac{1}{(2\pi i)^n}\int_{\gamma_0}\eta$, so the period is holomorphic.
    
\end{proof}  

\subsubsection{}
We have that the isomorphism \ref{eq: DRR for BCOV bundle} is realised as $DRR(\eta_{BCOV}^{12})=F(s)\eta_0^\frac{\chi}{m}$ for some holomorphic function $F(s)$. It is expected that the isomorphism $DRR$ lifts the arithmetic Grothendieck-Riemann-Roch of Gillet-Soul\'e, hence we would have a metric interpretation of this $F(s)$. For the purpose of this text we make the following definition: 
\begin{definition}\label{def:DRR function}
    In the context above we define the non-zero function $|DRR|:S\to \R$:
\begin{equation}\label{DRR-def}|DRR(s)|^2=\frac{||\eta_{BCOV}^{12}||^2_{Q,BCOV}}{||\eta_0^\frac{\chi}{m}||^2_{L^{\frac{2}{m}}}}=\frac{\left(\tau_{BCOV}\cdot ||\eta_{BCOV}||^2_{L^2,BCOV}\right)^{12}}{||\eta_0^\frac{\chi}{m}||^2_{L^{\frac{2}{m}}}}
\end{equation}
\end{definition}

Since we still lack the metric comparison in the functorial Riemann--Roch theorem, we are content to give the following lemma. We do not claim any originality on its proof as it is common knowledge which we recall for the convinience of the reader.  

\begin{lemma}\label{lema:logDRR is psh}
    
    Suppose that $\log \tau_{BCOV}$ on $(\D^\ast)^r \subseteq S$ is $L^2_{loc}$. Then $$|DRR(t)|=|\varphi(t)| \cdot \prod |t_i|^{2a_i},$$ where $\varphi(t)$ is an invertible holomorphic function on $\D^r$ and $a_i \in \R.$ If the family extends to a projective family over $\D^r$ we have that $a_i \in \Q.$
\end{lemma}

\begin{proof}
We write $u=\log|DRR|$. Since $-dd^c \log \|\eta_{BCOV}\|^{-2}_{Q,BCOV}$ computes the first Chern of the Quillen BCOV-bundle, and $-\frac{1}{m}dd^c \log \|\eta_0^{\frac{\chi}{m}}\|_{L^{\frac2m}}^{-2} $ computes the Weil--Petersson form, thanks to the holomorphic anomaly equation encapsulated in item \ref{item1:QuillenBCOV} of Prop.\ref{prop:QuillenBCOV}, we know that $u$ is pluriharmonic.

If we write $\omega=i(\partial u-\overline{\partial}u)$ we have that $d\omega=0$ and its de Rham cohomology class is characterized by the periods $a_j=\int_{\gamma_j}\omega$. This means that there exists a real-valued $g$ such that $dg +\sum_j \frac{a_j}{2\pi}d\log|t_j|=\omega$. This gives us that the function $u-\sum_j\frac{a_j}{2\pi}\log|t_j|-ig$ is holomorphic. So $u-\sum_j\frac{a_j} \log|t_j|$ is the real part of a holomorphic function $F$ on the punctured polydisk. 

The holomorphic function $F$ is also $L^2_{loc}$. Indeed, $|F|^2 = 2 \operatorname{Re}(F)^2 - \operatorname{Re}(F^2)$. The first one is $L^2_{loc}$ by the above, and the second one is $L^2_{loc}$ by Fubini and computing with the one-variable Laurent expansion of $F^2$. Hence by the $L^2$ Riemann removable singularity theorem, $F$ extends holomorphically to the  the full polydisk.

Finally, if the family can be completed to a projective family, by verifying against one-parameter degenerations we see that by \cite[Proposition 6.6]{DURF2} all $a_i \in \Q.$

\end{proof}

We remark that the condition that the BCOV invariant is $L^2_{loc}$ is automatically satisfied if the family extends to a projective family over a projective base. This is a consequence of the Schmid nilpotent orbit theorem and the arithmetic Riemann--Roch theorem, but we do not provide the details.

\subsubsection{} 
By construction the function $|DRR|$ is the $24$-th power of a  nonzero holomorphic function. With this in mind, we make the following definition. 

\begin{definition}
    When $|DRR|$ is the absolute value of a holomorphic function, it can be written as 
    $|C\exp(24(-1)^n F_1^B)|$
for a unique function $F_1^B$ which has the form 
    \begin{displaymath}
        F_1^B= \sum a_i \log t_i + \hbox{holo without constant term}. 
    \end{displaymath}
    We refer to $F_1^B$ as the genus one B-side potential. 
\end{definition}

\section{Mirror symmetry for a product}\label{sect: MS products}
In this section we study the behavior of the BCOV-invariant and the conjecture under products. This presents an interesting first step into studying how the set up of the refined $BCOV$ conjecture behaves under certain constructions that appear in algebraic geometry. These formulas relating the BCOV invariant of $X,Y$ to $X\times Y$ and their mirrors for the $A$-side potentials, allow for more examples where mirror symmetry is satisfied.  

The results of this section amount to the following theorem: 
\begin{theorem}\label{thm: M-S products}
    Let $X,Y$ be smooth projective varieties with torsion canonical bundle of dimensions $n,m$, then:
    \begin{itemize}
        \item $F_{1,X\times Y}^B=(-1)^n\chi(X)F_{1,Y}^B+(-1)^m\chi(Y)F_{1,X}^B$
        \item  $F_{1,X\times Y}^A=\chi(X)F_{1,Y}^A+\chi(Y)F_{1,X}^A$
    \end{itemize}
    where $F_1^B$ denotes the function computed with respect to the product family.
\end{theorem}
As the proof will show part $1$ of the above actually works more generally for compact Kähler manifolds with vanishing $c_1$. 
\subsection{The BCOV invariant of a product}
Unlike for quotients, the BCOV invariant has a particularly simple expression for products in terms of the individual BCOV invariants. For this subsection we recall the following notation; for a topological space $Z$ we denote by $\chi(Z)=\sum (-1)^k b_k$ (resp. $\chi'(Z)=\sum (-1)^k k b_k$) the topological (secondary) Euler characteristic. \\

\begin{proposition}\label{prop: BCOV invariant of product}
    Let $X,Y$ be compact Kähler manifolds with trivial Chern class, of dimension $n$ and $m$. Then 
    \[\tau_{BCOV}(X \times Y)= \tau_{BCOV}(X)^{\chi(Y)} \cdot \tau_{BCOV}(Y)^{\chi(X)}\cdot \left(\frac{1}{2\pi}\right)^{\left(n m\chi(X)\chi(Y)\right)}\]
\end{proposition}

For the proof, we assume for the moment that we endow $X,Y$ with a Ricci-flat Kähler metric, and we take the product metric on $X\times Y$. It is straightforward to verify that $A(X\times Y)=A(X)^{\chi(Y)}A(Y)^{\chi(X)}$, and we leave this to the reader. For the other factors, we begin by recording the following Lemma.\\
\begin{lemma}\label{lemma:Factor B for product}
We have 
\[
B(X \times Y) = B(X)^{\chi(Y)} \cdot B(Y)^{\chi(X)}\cdot\left(2\pi\right)^{\left(nm \chi(X)\chi(Y)\right)}.
\]
\end{lemma}

\begin{proof}    
    For a complex space $Z$ of complex dimension $N$, let $v_k=\log(\vol(H_k(Z,\Z),\omega_Z))$ and $b(Z)=\log B(Z)=\frac{-1}{2}\sum_{k=0}^{2N}(-1)^kkv_k$. With this definition we see that the definition of $B$ in \eqref{def:B-term} takes the form $b(Z) = \log B(Z).$
    We introduce the auxiliary polynomial 
    \begin{equation}\label{eq:auxpolyP}P_Z(t)=\frac{-1}{2}\sum_{k=0}^{2N}(-1)^kv_kt^k,
    \end{equation} so that $b(Z)=P'(1)$. Observe that, by item \ref{item3:renormalized L^2 norm} in \ref{prop:renormalized L^2 norm}, we have  $v_k + v_{2N-k}=- 2b_k N \log(2\pi)$ so that 
    \begin{equation}\label{eq:P(1)}P(1)=\frac{-1}{2}\sum (-1)^k v_k=\frac{\log(2\pi)}{2}\left( (-1)^N N b_N +N \sum_{k=0}^{N-1} (-1)^k (b_k + b_{2N-k})\right)=\frac{N}{2}\chi(Z)\log(2\pi).\end{equation}

Since we take the product metric, using that we have an isometry $$\det(A\otimes B)=\det(A)^{\rk(B)}\otimes \det(B)^{\rk(A)},$$ we find that $$\vol(H_{k_1}(X,\Z)\otimes H_{k_2}(Y,\Z),\omega_X\otimes \omega_Y))=\vol(H_{k_1}(X,\Z),\omega_X)^{b_{k_2}(Y)}\vol(H_{k_2}(Y,\Z),\omega_Y)^{b_{k_1}(X)},$$ 
where $b_k$ are the Betti numbers. Using this and the Künneth formula we deduce that 
\begin{equation}\label{eq:kunnethvol} v_k(X\times Y))=\sum_{k_1+k_2=k}b_{k_1}(X)v_{k_2}(Y)+v_{k_1}(X)b_{k_2}(Y)
\end{equation}

Likewise introduce $\chi_Z(t)=\sum (-1)^k b_k t^k.$  Then the expression \eqref{eq:kunnethvol} translates to 
\begin{equation}\label{eq:P-polynomial} 
P_{X\times Y}(t)=P_X(t)\chi_Y(t)+\chi_X(t)P_Y(t)
\end{equation} Recalling the expression \eqref{eq:P(1)}, the identities $\chi_Z(1)=\chi(Z)$ (topological Euler-Poincaré characteristic) and $\chi'_Z(1)=\chi'(Z)$ (secondary topological Euler-Poincaré characteristic), we can compute: 
$$b(X\times Y)=P_{X\times Y}'(1)=P_X'(1)\chi(Y)+P_Y'(1)\chi(X)+\frac{1}{2}(n \chi(X)\chi'(Y)+m \chi(Y)\chi'(X))\log(2\pi)$$

Lastly, Poincar\'e duality also shows that $\chi ' (X)=n\chi(X)$ so we get the desired result.
\end{proof}

We now proceed to find an analogous description of $\log(\mathcal{T}_{BCOV}(X\times Y))$, which together with the previous observations give proposition \ref{prop: BCOV invariant of product}.

\begin{lemma}\label{lemma:Factor T for product}
We have 
\[
\mathcal{T}_{BCOV}(X \times Y) = \mathcal{T}_{BCOV}(X)^{\chi(Y)} \cdot \mathcal{T}_{BCOV}(Y)^{\chi(X)}.
\]

\end{lemma}
\begin{proof} We follow the same strategy as in the previous Lemma. By \cite[Theorem 3.1]{ray1973analytic} we have $\sum_{p=0}^{N}(-1)^pt_p(Z)=0$ where $t_p(Z)=\log(T_p(Z))=\log(T(Z,\Omega^p))=\sum_q (-1)^q q\log \det(\Delta^{p,q})$. 
Moreover, the Künneth type formula of  \cite[Theorem 3.3]{ray1973analytic}  (taking the trivial characters)  gives us:
$$t_p(X\times Y)=\sum_{p_1+p_2=p}[t_{p_1}(X)\chi_Y(1,\Omega^{p_2})+\chi_X(1,\Omega^{p_1})t_{p_2}(Y)]. $$
Again using polynomial notation we introduce:
\begin{itemize}
    \item $Q_Z(t)=\sum_{p=0}^{N}(-1)^pt_p(Z)\cdot t^p$.
    \item $R_Z(t)=\sum_{p=0}^{N}(-1)^p\chi(Z,\Omega^p_Z)t^p$.
\end{itemize} 
Notice that $\log\mathcal{T}_{BCOV}(Z)=Q_Z'(1)$. The above results of Ray-Singer imply:
$$Q_{X\times Y}(t)=Q_X(t)R_Y(t)+Q_Y(t)R_X(t)$$ also $Q(1)=0$ and $R_Z(1)=\sum_p (-1)^p\chi(\Omega_Z^p)=\sum_{p,q}(-1)^{p+q}h_Z^{p,q}=\chi(Z)$. Thus we get 
$$\log(\mathcal{T}_{BCOV}(X\times Y))=\log(\mathcal{T}_{BCOV}(X))\chi(Y)+\log(\mathcal{T}_{BCOV}(Y))\chi(X)$$
\end{proof}

\subsubsection{Distinguished sections of a product} 

To compute $|DRR|$ of the product, we furthermore need to discuss how distinguished sections behave under products.

\begin{proposition}\label{prop:iso lines and sections of product}
There are two canonical isomorphisms which we express fiberwise as:

\begin{enumerate}
 \item \begin{equation}\label{eq:bcovlineproduct}\lambda_{BCOV}(X\times Y)=\lambda_{BCOV}(X)^{\chi(Y)} \otimes\lambda_{BCOV}(Y)^{\chi(X)}.
 \end{equation}
    \item \begin{equation}\label{eq:canonicalbundleprod}
    H^0(K^{m_{X\times Y}}_{X \times Y })^{\frac{\chi(X \times Y)}{12 m_{X \times Y}}}=\left(H^0(K^{m_X}_{X})^{\frac{\chi(X)}{12 m_{X}}}\right)^{\chi(Y)} \otimes \left(H^0(K^{m_Y}_{Y})^{\frac{\chi(Y)}{12 m_{Y}}}\right)^{\chi(X)}.
\end{equation}
Moreover, these are isometries for the $L^2$ or $L^{2/m}$-metrics. Under these isomorphisms choices of distinguished sections are related in the following way:
    \item For $\lambda_{BCOV}(X\times Y) $ we have 
    $$\eta_{BCOV}(X\times Y)\mapsto \eta_{BCOV}(X)^{\chi(Y)}\otimes \eta_{BCOV}(Y)^{\chi(X)}$$
    \item For $ H^0(K^{m_{X\times Y}}_{X \times Y })^{\frac{\chi(X \times Y)}{12 m_{X \times Y}}}$ we have 
    $$\eta_0(X\times Y)^\frac{\chi(X\times Y)}{12 m_{X\times Y}}\mapsto \left(\eta_0(X)^\frac{\chi(X)}{12 m_{X}}\right)^{\chi(Y)}\otimes \left(\eta_0(Y)^\frac{\chi(Y)}{12 m_{Y}}\right)^{\chi(X)}$$

\end{enumerate}
\end{proposition}
We remark that the distinguished sections are sometimes only defined up to a constant, and the above identification of said sections are up to a constant of the base. 

Now, the first point of Theorem \ref{thm: M-S products} is a consequence of the computation of the product of the BCOV invariant in Proposition \ref{prop: BCOV invariant of product}, the above Proposition \ref{prop:iso lines and sections of product} and the very definition of $F_1^B$ via $|DRR|$ via \ref{def:DRR function}. We now explain the proof of \ref{prop:iso lines and sections of product} in the following.

\subsubsection{} An argument as for what we did for BCOV-torsion but with coefficients as vector spaces (or bundles depending if we work over a fixed variety or over a family) yields the first isomorphism in Proposition \ref{prop:iso lines and sections of product}. 

In this case the polynomial we use is slightly different, as we only weight $\det H^{p,q}$ by $p$ in the definition of the BCOV line. Thus, by differentiating $P_{BCOV}(t)=\sum_p(-1)^p(\sum_q (-1)^q\det H^{p,q})t^p$, we find $\lambda_{BCOV}=\frac{d}{dt}P_{BCOV}(t)|_{t=1}$. It is an  $L^2$-isometry since the K\"unneth decomposition is. We are implicitly using an isometric trivialization $P(1)=(\C,|\cdot|)$ coming from Serre duality as in \cite[Prop. 3.4]{Fang-Lu-Yoshikawatorsionthreefolds}, if one chooses the trivialization coming from Poincar\'e duality some constant factors of $2\pi$ will appear (see \ref{lemma:Factor B for product}).

\subsubsection{} If we have two families $X/(\D^\ast)^{d_1}$, $Y/(\D^\ast)^{d_2}$ we can consider the product family $X\times Y$ over $(\D^\ast)^{d_1+d_2}$. In this setting the monodromy action is compatible with the Künneth decomposition so unipotency is preserved both for $X\times Y$ and their canonical coverings.

Moreover, under the K\"unneth decomposition any monodromy logarithm around a puncture is of the form $N=N_1\otimes 1+1\otimes N_2$ where $N_1,N_2$ are the monodromy logarithms of the starting families. This means that the Künneth decomposition is an isomorphism for the limit mixed Hodge structure. In particular, we have: 
\begin{equation}\label{eq:weightdecom}
    \Gr^W_lH^k(X\times Y)=\sum_{k_1+k_2=k}\sum_{l_1+l_2=l} \Gr^W_{l_1}H^{k_1}(X)\otimes \Gr^W_{l_2}H^{k_2}(Y)
\end{equation}
If the limit mixed Hodge structures associated to $X,Y$ are Hodge-Tate the right hand side of \eqref{eq:weightdecom} is a sum of Tate structures of type $\Z(-l_1)\otimes \Z(-l_2)=\Z(-l)$ (the integral structure of the graded pieces is preserved by K\"unneth). Hence the LMHS corresponding to the left hand side of \eqref{eq:weightdecom} is also Hodge-Tate. This description of the LMHS as a sum of tensor products of the LMHS of the factors gives that the appropriate sums of tensor products of distinguished sections give the distinguished sections of the family $X\times Y$ as in the third point of Proposition \ref{prop:iso lines and sections of product}.

\subsubsection{} We finally consider the last point of Proposition \ref{prop:iso lines and sections of product}. Denote by $m_X,m_Y$ the orders of $K_X$ and $K_Y$, and $m_{X \times Y}=\operatorname{lcm}(m,n)$. Since the canonical map of $\Q$-lines is an isometry $H^0(K_{Z}^{pm_Z})_{L^{2/{pm_Z}}}^\frac{1}{pm_Z}=H^0(K_Z^{m_Z})_{L^{\frac{2}{m_Z}}}^{\frac{1}{m_z}}$, there is a canonical isometry 
\begin{equation}
    H^0(K^{m_{X\times Y}}_{X \times Y })^{\frac{\chi(X \times Y)}{12 m_{X \times Y}}}=\left(H^0(K^{m_X}_{X})^{\frac{\chi(X)}{12 m_{X}}}\right)^{\chi(Y)} \otimes \left(H^0(K^{m_Y}_{Y})^{\frac{\chi(Y)}{12 m_{Y}}}\right)^{\chi(X)}.
\end{equation}
An argument as for the BCOV line shows that the product of distinguished sections is the distinguished section of the canonical bundle, now possibly up to a constant. The key here is again that, since the monodromy is unipotent, the limit mixed Hodge structures (and hence their splittings) behave well under products. Moreover, notice that from the aforementioned isometry $H^0(K_{Z}^{pm_Z})^\frac{1}{pm_Z}=H^0(K_Z^{m_Z})^{\frac{1}{m_z}}$ we get that to compute the norm of the section $\eta_0$ we can choose either the degree $m_{X\times Y}$ canonical cover or the degree $m_X\cdot m_Y$ cover given by the product of the canonical covers of $X$ and $Y$.

\subsection{The $F_1^A$ of a product}\label{subsect: GW products}

In this subsection we prove the behavior of  $F_1^A$ under products of smooth projective varieties with torsion canonical bundle. 
We will split the argument into two cases according to  wether the product of hodge numbers $h^{1,0}(X)h^{1,0}(Y)$ vanishes or not. The leading terms can be handed uniformly:

\subsubsection{}Recall that, in general, we denote by general term by
\begin{displaymath} N_{1,0}^{X}:=\frac{-1}{24}\int_{X}c_{n+m-1}(X)\cap 2\pi i\tau.
\end{displaymath}
For a product $X \times Y$ the Kähler parameter $\tau$ splits, by the K\"unneth decomposition, as $\tau=\tau_X+\tau_Y+m$ where $\tau_X,\tau_Y$ are Kähler parameters in $X,Y$ and $m$ is a ''mixed term'' i.e. a cohomology class coming from tensoring classes in $H^1(X)$ and $H^1(Y)$. 

Using the Whitney sum formula and the fact that the mixed terms capped with the chern classes are not of top degree we get the following additivity of the leading terms: 
\begin{equation*}
\begin{split}
    N_{1,0}^{X \times Y} & =\int_{X\times Y}c_{n+m-1}(X\times Y)\wedge (\tau_X+\tau_Y+m)\\&=
    \int_Yc_{m}(Y)\int_X c_{n-1}(X)\wedge \tau_X+\int_Xc_n(X)\int_Yc_{m-1}(Y)\wedge \tau_Y\\
    &=\chi(Y)N_{1,0}^X+\chi(X)N_{1,0}^Y
    \end{split}
\end{equation*}

\subsubsection{} We need to prove the analogue of the above additivity in general. Let us begin by assuming that either $h^{1,0}(X)$ or $h^{1,0}(Y)$ vanishes. In this case we have Beherend's formula (\cite{BeherendProduct}) giving the following equality of Gromov--Witten classes:

\begin{equation}\label{eq: beherendprod}
    I_{g,n,\beta}(\gamma_X\otimes \gamma_Y)=I_{g,n,\beta_X}(\gamma_X)\cup I_{g,n,\beta_Y}(\gamma_Y) 
\end{equation}
This equation follows from the existence of the following diagram with cartesian squares :

\[\begin{tikzcd}
	{M_{g,n}(X\times Y)} & P & {M_{g,n}(X)\times M_{g,n}(Y)} \\
	{D_{g,n}} & Q & {\mathfrak{M}_{g,n}\times \mathfrak{M}_{g,n}} \\
	& {M_{g,n}} & {M_{g,n}\times M_{g,n}}
	\arrow["h",from=1-1, to=1-2]
	\arrow[from=1-1, to=2-1]
	\arrow[from=1-2, to=1-3]
	\arrow[from=1-2, to=2-2]
	\arrow[from=1-3, to=2-3]
	\arrow[from=2-1, to=2-2]
	\arrow[from=2-2, to=2-3]
	\arrow[from=2-2, to=3-2]
	\arrow[from=2-3, to=3-3]
	\arrow["\Delta",from=3-2, to=3-3]
\end{tikzcd}\]
This implies, see \cite{LogProductFormula}, the equality between virtual fundamental classes:
$$h_\ast [\amalg_\beta M_{g,n}(X\times Y)]^{virt}=\Delta^!([M_{g,n}(X)]\times [M_{g,n}(Y)]^{virt}) $$
In general this is for moduli spaces of stable curves without fixing a homology class. When passing to homology one needs to remember that the degree $2$ homology of $X\times Y$ is not the product of the degree $2$ homology of $X$ and $Y$ unless one of these factors has vanishing first Betti number, as is the case for us. 

\subsubsection{} We apply \eqref{eq: beherendprod} to $\alpha=\alpha_X\otimes 1_Y+1_X\otimes\alpha_Y\in H^2(X)\oplus H^2(Y) =  H^2(X\times Y)$, and obtain
$$\langle \alpha\rangle_{1,1,\beta}=I^Y_{1,1,\beta_Y}(1_Y)\cdot \langle \alpha_X\rangle _{1,1,\beta_X}+I^X_{1,1,\beta_X}(1_X)\langle\alpha_Y\rangle_{1,1,\beta_Y}$$
Here we use $1_X$ to denote the class corresponding to the generator of $H^0(X,\Z)$.

When $\beta_X\neq 0$ the classes $I_{1,1,\beta_X}(1_X)$ are zero by a standard application of the reduction axiom (see \cite[Subsubsect 7.3.1]{KatzCoxMirrorBOOK}). This gives that the only GW invariants that we get are of the form  $\langle \alpha \rangle^{X\times Y}_{1,1,\beta}$ such that either $\beta_X=0$ or $\beta_Y=0$.

For simplicity we assume $\beta_Y=0$. By the arguments above we want to compute $\langle \alpha\rangle^{X\times Y}_{1,1,\beta=\beta_X}=I^Y_{1,1,0}(1_Y)\langle \alpha \rangle^X_{1,1,\beta_X}$. By the  divisor axiom for Gromov--Witten invariants  (see \cite[Subsubsect 7.3.1]{KatzCoxMirrorBOOK}) we obtain:
$$ \left(\int_{\beta}\alpha\right) N^{X\times Y}_{1,\beta}=I^Y_{1,1,0}(1_Y)\left(\int_{\beta_X}\alpha_X\right)N^X_{1,\beta_X}$$ 
and we can divide away $\int_\beta \alpha = \int_{\beta_X} \alpha_X$ by choosing it to be non-zero. 

To finish, we know 
$$I^Y_{1,1,0}(1_Y)=\int_Yc_n(T_Y) 1_{M_{1,1}}$$  
by applying the definition of this specific GW class (see \cite[Subsubsect 2.4.5]{kontsevich1994gromov}).

Thus, we conclude that
$$N^{X\times Y}_{1,\beta_X}=\chi(Y)N^X_{1,\beta_X}(X)$$
which entails the additivity for $F_1^A$ (i.e. part 2 of Theorem \ref{thm: M-S products}).

\subsubsection{} We now suppose that both $h^{1,0}(X)$ and $h^{1,0}(Y)$ are non-zero. We remark that for any $Z$ with $c_1(Z)=0$, by the Beauville--Bogomolov decomposition (see \cite{BeavilleVarietesClasschernnulle}, and also the formulation of Kawamata \cite[Theorem 8.3]{kawamata85}), there exists a cartesian diagram:
\begin{equation}\label{eq:beauville-bogomolov}
    \xymatrix{
    Z' \times B \ar[r] \ar[d]^{q} & B \ar[d] \\
    Z \ar[r] & A}
\end{equation}
where $A$ is the Albanese and $B \to A$ is a finite abelian (e\'tale) cover. Since $\dim(B)=\dim(A)=h^{1,0}(Z)$ and the Euler characteristic of an abelian variety vanishes, we get $\chi(Z)=0$ whenever $h^{1,0}(Z)$.

This leads us also to prove the following stronger version of the additivity of $F_1^A$ for $h^{1,0}(X)$ and $ h^{1,0}(Y)$ both $ \geq 1$, which shows that both sides are zero: 
\begin{proposition}\label{prop: GW are zero for b_1 nonzero}
    Suppose that $Z$ is a projective manifold with $c_1(Z)=0$ and $h^{1,0}(Z) \geq 2$. Then $\chi(Z)=0$ and $F_{1,Z}^{A}=0$. 
\end{proposition}
\begin{proof}
    We use the presentation in \eqref{eq:beauville-bogomolov}. By deforming $B$ and $A$, we can suppose that $A$ is simple, i.e. not isogenous to a product of abelian varieties. Simplicity implies the image of any elliptic curve in $A$ is constant. In particular it also follows that any 1-pointed nodal genus one map $(C,p) \to Z$ has constant image in $A$, i.e. lies in the fiber of the Albanese morphism $Z \to A$.

    Given a map $(C,p)\to Z$ with image $a \in A$ and picking a lift $b\in B$ of $a$, is the same thing as a map $(C,p) \to Z \times_A B = Z' \times B$ which is constant on $B$. It gives us an isomorphism of moduli stacks:
    $$M_{1,1}(Z)\times_A B=M_{1,1}(Z')\times B.$$
    
    Restricting this to effective curve classes $\beta$ in $H_2(Z)$, this decomposition takes the form 
    $$M_{1,1}(Z, \beta)\times_A B=\bigcup_{\beta'}M_{1,1}(Z', \beta')\times B,$$
    where $\beta'$ runs over the effective curve classes in $H_2(Z')$ whose image in $H_2(Z)$ is $\beta$. Notice that since the left hand side is a proper finite type DM stack the union on the right hand side is finite. Moreover by construction of this decomposition their perfect obstruction theories are compatible. By \cite[Prop. 5.10]{FantechiBeherendIntrinsic} we have that 
    \begin{displaymath}
        q^\ast [M_{1,1}(Z, \beta)]^{virt} = \sum_{\beta'} [M_{1,1}(Z', \beta') \times B]^{virt}
    \end{displaymath}
    Since $q_\ast q^\ast = \deg[B/A]$, it is enough to prove that the classes $[M_{1,1}(Z',\beta')\times B]^{virt}$ vanish, and then apply the divisor axiom to show that $N_1^{X}$ vanishes. We conclude by the following equalities. 

    \begin{equation*}
    \begin{split}
    [M_{1,1}(Z',\beta')\times B]^{virt} &=
    [M_{1,1}(Z', \beta') \times_{M_{1,1}} M_{1,1}(B, 0)]^{virt}\\
    &=\Delta^![M_{1,1}(Z',\beta')\times M_{1,1}(B,0)]^{virt}\\ &=\Delta^! \left( [M_{1,1}(Z',\beta')]^{virt}\times [M_{1,1}(B,0)]^{virt}\right)
    \\ &=\Delta^!\left([M_{1,1}(Z', \beta')]^{virt} \times\left( c_{\text{top}}(R^1p_\ast \mu^\ast T_B ) \cap [B]\right) \right).
    \end{split}
    \end{equation*}
    where the first equality follows from the description of the moduli stacks , the second from \cite[Prop. 5.10]{FantechiBeherendIntrinsic}, the third is \cite[Prop. 5.7]{FantechiBeherendIntrinsic} and the fourth is \cite[Prop 5.6]{FantechiBeherendIntrinsic}. We can finish thanks to the explicit description of the perfect obstruction theory for the spaces of stable curves (here $p:\mathcal{C}\to M_{1,1}(B,0) = M_{1,1} \times B$ is the universal curve and $\mu:\mathcal{C}\to B$ its universal evaluation map).
    
    Since $T_B$ is trivial, the obstruction bundle $R^1p_\ast \mu^\ast T_B= \left(R^1 p_\ast \Ocal_{\mathcal{C}}\right)^{\dim B}$. The bundle $\lambda = R^1 p_\ast \Ocal_{\mathcal{C}}$ is the dual of the Hodge bundle on $M_{1,1}$, and we have that $c_{\text{top}}(R^1p_\ast \mu^\ast T_B)= c_1(-\lambda)^{\dim B}$. This is zero if $\dim B \geq 2$ because $M_{1,1}$ is one dimensional. This is the case since $\dim B = h^{1,0}(Z)$. 
\end{proof}

\section{Lattice polarized K3 surfaces and maximally unipotent points}\label{section:latticepolarizedk3}

For Enriques surfaces, a mirror symmetry statement and a computation of the BCOV invariant have been carried out in the literature. The generating series of Gromov--Witten invariants, recalled below, has been computed by Oberdieck \cite{Oberdieck2025KMformula} in terms of a certain Borcherds lift $\Phi$. This completed a conjecture of Klemm--Mariño \cite{MR2391190}. 

The corresponding analytic torsion was computed by Yoshikawa in \cite[Theorem 8.3]{YoshikawaK3inventiones}. It follows from this, combined with \cite[Theorem 9.2, Eq 9.2]{DaiYoshikawaLogEnriques}, that for an Enriques surfaces $Y$ we have :
\begin{equation}\label{eq:yoshikawaBCOV}
    \tau_{BCOV}(Y)= C\cdot \|\Phi \|_{\operatorname{Pet}}^{\frac12}
\end{equation}
where $\|\Phi\|_{\operatorname{Pet}}$ is the Petersson norm of the same Borcherds lift as in the Klemm--Mariño formula and $C$ is independent of $Y$.

As we can see both computations are related to some power of the Borcherds modular form. Below, we recast Yoshikawa's result in the formalism of distinguished sections and show that these results put together prove $F_1^A=-F_1^B$ for Enriques surfaces as mentioned in the introduction.  More precisely our main result is the following: 
\begin{theorem}\label{thm:mirror theorem Enriques}
    There exists a strongly unipotent family of Enriques surfaces $f:X\to (\mathbb{D}^\ast)^{10}$ where the canonical section of $f_\ast(K_{X/S}^2)$ is given by the nilpotent orbit associated to the family.  
    
    Moreover, the family admits a change of coordinates defined in terms of holomorphic periods
    $$q:(\D^\ast)^{10}\to (\D^\ast)^{10},$$
    such that $F_1^B(q(s))=-F_1^A(\tau)$. 
\end{theorem}

\subsection{Moduli of $M$-polarized K3}\label{subsect: moduli polarized K3}

The study of the Enriques moduli can be done in general for so-called lattice polarized $K3$ surfaces and coincides with an explicit period domain. In the present section we study the tube domain presentation of said moduli space.

\subsubsection{}We follow \cite{Dolgachev1996, Nikulin1980, YoshikawaK3inventiones, AlexeevEngel2025LatticePolarizedK3} and the summary of \cite[Ch. 8]{Compactcomplexsurfaces}. Let $X$ be an algebraic K3 surface. Recall that $H^2(X, \Z)$ is isomorphic to the K3-lattice  
\[\Lambda_{K3}= U^{\oplus 3} \oplus E_8(-1)^{\oplus 2}.
\]
Here $U$ is the even unimodular hyperbolic lattice, generated by $e,f$ and $(e,e)=(f,f)=0, (e,f)=1$ and $E_8$ is an even positive definite lattice of rank 8. The K3-lattice is of signature $(3,19)$. 

Recall that we are choosing a renormalization of the $L^2$ inner product to make it so that it induces a polarized Hodge structures as in \cite[Def. 2.9 bis]{MixedHodgestructures}. This means that our bilinear form in integral cohomology is $\frac{1}{(2\pi i)^2}$ times cup product. Thus, in our context, the isomorphism $(H^2(X,\Z),(,)_{L^2})\cong \Lambda_{K3}$ has an implicit factor of $(2\pi i)^2$. An alternative way of thinking about this is that, with our renormalization we study how $H^2(X,\Z(1))$ varies, then the pairing takes values in $\Z(0)$.

For better exposition we present the results on the moduli space without this normalization in the lattice. The above mentioned technicality only becomes important when we compute explicitly certain norms of distinguished sections and when we state the equalities in \ref{eq: graded pieces of K3 LMHS} between graded objects and pure Tate structures.

\subsubsection{} Suppose $M$ is a lattice of rank $r$ of signature $(1,r-1)$, and there is a primitive embedding $M \to \Lambda_{K3}$. Associated to $M$ is $M^\perp$, a lattice of signature $(2,20-r)$. The period domain of $M$ is a symmetric domain of type IV:
\[\mathcal{D}_M= \{[w] \in \P(M^\perp_\C) ,\; \quad (w,w) =0, \quad (w,\overline{w})>0\}=\mathcal{D}_L^+ \sqcup \mathcal{D}_L^-.
\]

Both of the two components lie in the compact dual $\mathcal{D}_L^\vee \subseteq \P(L_\C)$ defined by the single quadric equation $(w,w)=0$. 

An element $[w] \in \mathcal{D}_M$ defines a Hodge structure of weight 2 and $h^{0,2}=h^{2,0}=1$ on $L$ by $F^2 \subseteq F^1 \subseteq F^0$ being given by $[\omega] \subseteq [w]^\perp \subseteq M^\perp_\C$. The first condition simply gives the inclusion $F^2 \subseteq F^1$ and the second condition is {equivalent} to that $\overline{H^{0,2}}=H^{2,0}$ and that the pairing polarizes the Hodge structure. Viewing $M$ as a trivial Hodge structure of pure type $(1,1)$, the decomposition $M_\C \oplus M^\perp_\C = \Lambda_{K3, \C}$ equips $\Lambda_{K3}$ with the structure of a polarized Hodge structure. 

\subsubsection{} Following \cite{AlexeevEngel2025LatticePolarizedK3}, we define a small cone $\mathfrak s \subseteq M_\R$ to be a connected component of the complement in 
\begin{equation}\label{eq:positivecone}
    C_{M}= \{y \in M_\R, (y,y) > 0 \} = C_{M}^+ \sqcup C_{M}^-
\end{equation}
of hyperplanes $\beta^\perp = \{x \in M_\R, (x,\beta)=0\}, \beta \in \Lambda_{K3} \setminus M^\perp$, such that $\beta^{2}=-2$ and $M+ \Z \beta$ is still Lorentzian. A $K3$-surface $X$ is $(M, \mathfrak s)$-quasi-polarized, or quasi-polarized by $(M, \mathfrak s)$, if there is a primitive lattice embedding $M \to \Pic(X) \subseteq H^2(X, \Z)$, such that  $\mathfrak s \subseteq M_\R$ is contained in the closure $\overline{K}_X$ of the K\"ahler cone $K_X$. We say it is $(M, \mathfrak s)$-polarized if $\mathfrak s \cap K_X \neq \emptyset$. The K3 surface is marked if there is an isometry $\phi: H^2(X, \Z) \simeq \Lambda_{K3}$ that takes $M$ to the fixed $M \subseteq \Lambda_{K3}$. The period of a marked quasi-polarized K3 surface determines a point in $\mathcal{D}_M$. 

There is a fine moduli space of marked $(M, \mathfrak s)$-polarized K3 surfaces, $\mathcal{M}_{M, \mathfrak{s}}^\circ$ constructed in \cite[Corollary 4.8]{AlexeevEngel2025LatticePolarizedK3}. It is a version of the case without $M$ considered in \cite{BurnsRapoport75}, and a technical correction of the first results of this nature in \cite[Corollary 3.2]{Dolgachev1996}. For this, consider the roots $\Delta(M^\perp) =\{\delta \in M^\perp, (\delta, \delta) = -2 \}$ and the corresponding wall $\delta^\perp = \{x \in M^\perp_\C, (x, \delta) =0\}.$ Then the period map 
\begin{equation}\label{eq:modulispaceMample-marked}\mathcal{M}_{M, \mathfrak s}^\circ \to \mathcal{D}_M^\circ = \mathcal{D}_M \setminus \bigcup_{\delta \in \Delta(M^\perp)} \delta^\perp
\end{equation}
is an isomorphism and we identify the two spaces. 

\subsubsection{Tube domains around 0-cusps} We now suppose that $\rk M \leq 19$, and furthermore suppose that $L=M^\perp=U \oplus K$ for some lattice $K$, necessarily of signature $(1,19-r)$.  Then, by \cite[Theorem 4.2, Remark 4.5]{Dolgachev1996} , if $C_K, C_K^\pm$ is defined as in \eqref{eq:positivecone} the map \[K_\R + i C_K \to \mathcal{D}_M\] 
given by $z= x+iy \mapsto \omega(z) = f + z - \frac{(z,z)}{2} e$ is an isomorphism and provides a tube domain description of the period domain. Here $f\in U$ is isotropic and $e\in U$ is such $(f,e)=1$. Observe that, by construction, if $k\in K$ we have 
\begin{displaymath}
    (k,z)=(k,\omega(z))
\end{displaymath}
since the decomposition $M^\perp=U\oplus K$ is orthogonal.

We recast \eqref{eq:modulispaceMample-marked} as a natural identification \[\mathcal{M}_{M, \sigma}^\circ= (K_\R + i C_K) \setminus \bigcup_{\delta \in \Delta(M^\perp)} \Phi^{-1} \delta^\perp\] 

\subsubsection{}\label{subsubseq:thewalls} If we write $\delta = ae+k+bf$, the condition that $z=x+iy$ is orthogonal to $\delta$ is, by taking real and imaginary parts, translated into \[a+(x,k) - b\frac{(x,x)}{2}+b \frac{(y,y)}{2}=0, \quad (y,k-bx)=0.\]
A manipulation shows that if $b \neq 0$, these equalities, together with $(\delta, \delta)=2ab  + (k,k)=-2$, becomes that $$(k-bx, k-bx) = (k,k)+b^2(y,y)-2b(k,x)=-2-2ab+b (y,y)+2ab=b^2(y,y)-2, \quad k-bx\in y^\perp.$$ For signature reasons, since $(y,y)>0$ the second condition implies that $b^2(y,y) < 2$. Then, if $b \neq 0$, $(y,y) < \frac{2}{b^2}\leq 2$, showing that if $(y,y)>2$ it will automatically avoid this wall. If $b=0$, the conditions on the tubular domain coming from removing $\delta^\perp$ are linear conditions defined by the linear walls $(y,k)=0$ and $(x,k)=-a$. We call  any connected component of the complement of $C_K$ against subspaces in $K_\R$ orthogonal to $-2$-vectors in $K$ a Weyl chamber.

\subsection{A mirror family of $M$-polarized $K3$}

Using the tube coordinates we construct a family of lattice polarized $K3$ surfaces and study the associated nilpotent orbit. This orbit is shown to converge to a Hodge--Tate degeneration around the $0$-cusps of the period domain which allows us to construct distinguished sections and a mirror map. In order to proceed we fix the component $C_K^+$ of $C_K$ and hence also the corresponding $\Dcal_M^+$.

\subsubsection{} We consider $\omega(0)$ lying in the compact dual, which has a Hodge filtration $F^\bullet(f): F^2 \subseteq F^1 \subseteq F^0$ given by $\C f \subseteq f^\perp \subseteq M^\perp_\C,$ also  consider a filtration $W_0\subseteq W_2 \subseteq W_4$ given by $\Z \cdot e \subseteq e^\perp \subseteq M^\perp$. 
Then 

\begin{equation}\label{eq:referenceHodgeTate}(F^\bullet(f), W_{\bullet, \Q}, M^\perp)
\end{equation} 
defines a mixed Hodge structure, with induced integral structures on the graded pieces:
\begin{equation}\label{eq: graded pieces of K3 LMHS}
    \Gr_0^W = \Z(0)= \Z \cdot e, \quad \Gr_2^W = \Z(-1)^{\oplus 20-r} = K, \quad \Gr_4^W = \Z(-2)= \Z \cdot f, 
\end{equation}

and it is hence Hodge--Tate.

\subsubsection{} For $k \in K$, define operators $N_k$ on $M^\perp_\C$ by 
\begin{equation}\label{eq:Nilpotentoperators}N_k(f)=k, \quad N_k(k')=-(k,k') e, \quad N_k(e)=0.
\end{equation}This is a nilpotent operator, so that $T_k = \exp(N_k)$ is unipotent, preserves the weight filtration $W_\bullet$ and $T_k(f) = f + k - \frac{(k,k)}{2} e$. Since the lattice is even, the last element is in fact in the lattice and $T_k$ preserves $M^\perp$.  The group generated by these $T_k$ are in fact the (integral) unipotent radical of the parabolic $P \subseteq O(M^{\perp})^+$ stabilizing the isotropic line $\Z \cdot e$. The corresponding Lie algebra is commutative, and isomorphic to $K$ via the map $k \mapsto N_k$. This hence determines a representation $T: K \to O(M^\perp)^+$.

\subsubsection{}
As operators on the periods of the Hodge filtrations, for $k \in K$, the action is given by \[
T_k [\omega(z)] = [\omega(z+k)]
\]
and we can consider the quotient of the VHS underlying the trivial local system $\Lambda_{K3}$ on $K_\R+iC_K^+$ by the lattice $K$. This hence provides us with a variation of Hodge structures on \[(K_\R+i C_K^+)/K  \subseteq K_\C/K\] 
The quotient $K_\C/K$ is naturally isomorphic to the torus \[ T_K := \Hom(K^\vee, \C^\ast) = K \otimes \C^\ast, \] where the isomorphism is given by $z \mapsto (m \mapsto \exp(2\pi i \langle m, k )\rangle),$ where $\langle -, -\rangle: K^\vee \times K \to \Z$ is the natural pairing.

\subsubsection{} We proceed to show that we can similarly descend the family. Consider a sector $K_\R + i W$, where $W$ is a Weyl chamber. The walls in $C_K$ from \ref{subsubseq:thewalls}, given by $(y, k)=0$ from $\bigcup_{\delta \in \Delta(M^\perp)} \delta^\perp$, are preserved by translation by $K$ and hence $K$ acts on $ K_\R + i W$ by translation. The points $K_\R + i W_{>2}$ where $(y,y)>2$ hence admit a family of K3 surfaces by the discussion in \ref{subsubseq:thewalls}. 

 The space $K_\R + i W_{>2}$ is simply connected analytic space, and by the pointwise global Torelli theorem, and the Main Lemma of \cite{BurnsRapoport75}, the action on the local system extends to an action of the whole family of K3 surfaces over $K_\R + i W_{>2}$. This lifted action satisfies the same cocycle condition as that of the action of $K$ on $\Lambda_{K3}$, and hence provides descent datum to descend the family to $(K_\R + i W_{>2})/K$.

 \subsubsection{} Pick a rational cone $\sigma \subseteq W$. By possibly shrinking it, we suppose it is a regular cone, i.e. we can choose an integral basis $k_1, \ldots, k_{20-r}$ of $K$, such that the new cone given by the convex hull $$\sigma = \operatorname{conv}(k_1, \ldots, k_{20-r}).$$ Then $K_\R + i \sigma$ is preserved by the action of $K$. Denote by $N_j = N_{k_j}$ the nilpotent operator corresponding to $k_j$. Then, for $z\in K_\C \simeq \C^{20-r}$, defining $N_z = \sum z_j N_j$, we have
\begin{equation}\label{eq:Nilpotentorbit}
\omega(z)=\exp(N_z) f = \exp\left(\sum z_j N_j\right) f.
\end{equation}  Hence  $\omega(z)$ is a nilpotent orbit with monodromy cone given by 
$$\sigma \simeq \left\{ \sum_{\lambda_j > 0} \lambda_j N_j  \right \}.$$
We recall this means that that if $\operatorname{\im} z_j$ is big enough $\omega(z)$ is a (pure) polarized Hodge structure,  and that the operators $N_j$ are commuting nilpotent operators satisfying Griffiths transversality.  If we pick the imaginary part such that $(y,y) > 2$, by \ref{subsubseq:thewalls} the corresponding part of the cone also lies in $W_{>2}$. The expression \eqref{eq:Nilpotentorbit} also shows that the tube domain itself is a type of nilpotent orbit.

\subsubsection{}\label{subsubsec:regularcone}

The regular cone $\sigma$ defines an isomorphism $K \simeq \Z^{20-r}$ and a toric partial compactification \begin{equation}\label{eq:K-torus} T_K = K_\C/K = \Hom(K^\vee, \C^\ast)= K \otimes \C^\ast \subseteq T_{K, \sigma}.
\end{equation}The monodromy determined by $K_\R + i \sigma \ni z \mapsto z+k$ along a ray generated by a primitive vector $k \in K$ is by construction $N_k$, and hence the monodromy around the $j$-th copy of $\C$ is $T_j=\exp(N_j)$. 

The ring of regular functions is $\C[\sigma^\vee \cap K^\vee],$ so that regular functions can be written as
\begin{displaymath}
    \sum_{m \in K^\vee, \langle m, \sigma \rangle \geq 0} c_m \chi^{m}
\end{displaymath}
where $\sigma^\vee \cap K^\vee$ is the set of all $m \in K^\vee$ such that $\langle m, \sigma \rangle \geq 0$. The expression $\chi^m$ is the natural map $\chi^m : K \otimes \C^\ast \to \C^\ast$. If we fix a basis $K \simeq \Z^{20-r}$, which is also the generators of the cone $\sigma$ and take the standard coordinates $t_i$ on $(\C^\ast)^{20-r}$ then $\chi^m(t)=t_1^{m_1}\ldots t_{20-r}^{m_{20-r}}$. The condition $m \cdot \sigma \geq 0$ is just that each $m_i \geq 0.$ By possibly having to rescale the torus, we think of $K_\R + i \sigma/K$ in $T_{K, \sigma}$ as a disc $\D^{20-r} \subseteq \C^{20-r}$.

\subsubsection{} The limit mixed Hodge structure of the constructed family $X \to \D^{20-r}$ associated to $\sigma$ is, by definition, determined by the Hodge flag corresponding to
\[
\lim_{\operatorname{\im} z \to \infty} e^{-\sum z_j N_j} \omega(z).
\]
But because of \eqref{eq:Nilpotentorbit} we see that $e^{-\sum z_j N_j} \omega(z)=F^\bullet(f)$. For any $N = \sum \lambda_j N_j, \lambda_j > 0$, we have by \eqref{eq:Nilpotentoperators} that \[N^2 f = -(\sum \lambda_j k_j, \sum \lambda_j k_j) e \neq 0\] since $(\sum \lambda_j k_j, \sum \lambda_j k_j)>0$, so $N^2 \neq 0$. This shows that the limit mixed Hodge structure is in fact \eqref{eq:referenceHodgeTate}, which  was Hodge--Tate by construction. We have proven the following proposition:\\

\begin{proposition}\label{prop:HodgeTatedegeneration}
    Suppose $M\subseteq \Lambda_{K3}$ is a primitive lattice of signature $(1,r-1)$, and that we are given a small cone $\mathfrak{s} \subseteq M_\R$. Suppose that $M^\perp = U \oplus K.$ Then there is a family of $(M, \mathfrak{s})$-polarized K3-surfaces $X \to (\D^\ast)^{20-r}$, which is of maximal variation and the limit mixed Hodge structure is Hodge--Tate. In particular, it is strongly unipotent and hence maximally unipotent.
\end{proposition}
The construction of the family a priori depends on the choice of a regular cone $\sigma \subseteq C_K$, but the conclusion does not. It would hence also make sense to talk about the whole Weyl chamber being Hodge--Tate.
\subsection{Distinguished section and mirror coordinates}

For a family as in Proposition \ref{prop:HodgeTatedegeneration}, the Hodge--Tate
limit mixed Hodge structure gives a distinguished section of \(F^n\), see
definition \ref{def:Hodge Tate distinguished sections}. 
We recall that these are sections  $(M^\perp)_\C^{p,2-p}$ corresponding to integral elements of $\Gr_{2p}^W$ via the isomorphisms $(M^\perp)_\C^{p,2-p} \simeq F^p \cap W_{2p} \simeq  \Gr_W^{2p} (M^\perp)_\C$. In this subsection we relate these constructions to the formulations used in the theory of orthogonal modular forms.

\subsubsection{} In the present polarized $K3$ situation, for $p=2$ this means the class $\eta \in F^2$ projecting onto $f$ in $\Gr_4^W$. Given a class $\eta$, one can obtain this by considering
\begin{equation}\label{eq:distF2}
    \frac{\eta}{(\eta,e)}.
\end{equation}
In tube coordinates we have $\omega(z) \in F^2$ and $(\omega(z),e)=1$
 so, \(\omega(z)\) is already the distinguished section. Notice that since 
 \begin{displaymath} W_{-4}H_2=\langle \gamma_0\rangle\cong\left(H^2/W_3H^2\right)^\vee =\langle \gamma_0\rangle\cong\left(H^2/W_2 H^2\right)^\vee =\Hom(\C \cdot f,\C),
 \end{displaymath}we can think of the generator $\gamma_0$ as the cycle that pairs as 1 with $f$ and zero with all other classes. Hence $e$ is the Poincaré dual to the homology cycle $\gamma_0$ so that $\int_{\gamma_0} \eta = (\eta, e)$. This shows that \eqref{eq:distF2} can likewise be understood as the, possibly more, familiar quotient $$\frac{\eta}{\int_{\gamma_0} \eta}.$$

The definition of $\gamma_0$ shows that it is an invariant cycle, since $N(W_{-4})=0$ for degree reasons.

In a similar way, the distinguished sections for $p=1$ correspond to elements $K \subseteq H^{1,1} \cap W_{2}=K_{\C}$ and for $p=0$ correspond to $e \in W_0.$

\subsubsection{}  We now describe the mirror coordinates. The extension
\[
    0\longrightarrow W_2/W_0\longrightarrow W_4/W_0\longrightarrow W_4/W_2\longrightarrow 0
\]
describes $W_4/W_0$ as an extension of $\Z(-2)$ by $K$. Since the quotients $\Gr_{2p}^W$ are trivial local systems, we identify them with their global sections and think of the above as an extension of variations of mixed Hodge structures. In general these are understood as the Jacobian of the  $\Hom$ mixed Hodge structure, cf.  \cite{CarlsonHainExtsensions} or \cite[Theorem 3.31]{MixedHodgestructures}.

In particular, we have the identity
\begin{equation}\label{eq: expression of ext MHS}
    \Ext^1_{\mathrm{MHS}}(\mathbb Z(-2),\mathbb Z(-1))=\Ext^1_{\mathrm{MHS}}(\mathbb Z(-1),\mathbb Z(0))=
    \mathbb C^\ast,
\end{equation}
this is explicitly described in \cite{DeligneLocalBehav} as being determined by an element $q\in \C^\ast$ such that the extension $H$  is generated by $e_1$ and $e_0+\frac{\log q}{2\pi i}e_1$ with $e_0$ generating $\Z(0)$ and $e_1$ generating $\Z(-1)$.

We can identify $W_4/W_0$ with the extensions of variation of mixed Hodge structures, 
\[
    E \in \Ext^1_{VHMS}(\Z(-2), K)=\Hom(\Z,K) \otimes \Ocal(\D^{20-r})^\ast = \operatorname{Hom}(K^\vee,\Ocal(\D^{20-r})^\ast).
\]
Hence, any $m \in K^\vee$ determines  $E_m:\D^{20-r} \to \C^\ast.$ These are described as follows. Take a lift of $f$ living in $F^2$ and write it as $f+z, z \in K_\C$. By \cite[Theorem 3.31]{MixedHodgestructures}, the extension is determined by the class of $z$. 

\subsubsection{}The image of an element $m\in K^\vee$ corresponds to the function $T_K = K_\C/K \to \C^\ast$ given by \begin{displaymath} E_m(z) = \exp(2\pi i \langle m, z\rangle).
\end{displaymath}
In the above expression we think of $z$ as an equivalence class in $K_\R + i \sigma/K\subset T_{K,\sigma}$.
 This map is nothing but, up to sign determined by the trivialization $\Gr_{4}^W M^{\perp} \simeq \Z$, the character $\chi^m$ of $T_K$ discussed in \ref{subsubsec:regularcone}. This can be seen by choosing a basis of $K$, using the fact that $K$ is the direct sum of pure Tate structures, and then using Deligne's description of \ref{eq: expression of ext MHS}.

 Recall also that the dual of the unique invariant cycle  $\gamma_0$ of $\Gr_{-4}H_2$ , is $e$. Then we find that 
\[
\langle m, z \rangle= \langle m, \omega(z) \rangle=\int_m  \omega(z)=\frac{\int_m\omega(z)}{\int_{\gamma_0}\omega(z)}.
\]
This is the usual representation of mirror coordinates as a quotient of particular periods.

\subsubsection{Orthogonal modular forms.}

We next recall the relation between two equivalent ways of viewing
orthogonal modular forms. On the one hand, the period domain carries the
tautological Hodge line bundle
\[
       \mathcal{L}:=\mathcal{F}^2 ,
\]
whose fiber over a period point $[\omega(z)]\in \mathcal D_M$ is the line
$\mathbb C \cdot \omega(z)$. This choice of $\omega(z)$ actually trivializes the bundle. 

For our purposes, a modular form of weight $k$
with respect to an arithmetic group $\Gamma$ is a $\Gamma$-invariant
holomorphic section of $\mathcal L^{\otimes k}$, possibly with a character. This is nothing but an expression of the form $F(z) \omega(z)^k$, invariant under $\Gamma.$ 

On the other hand, in the classical analytic description one writes such
a section as a holomorphic function satisfying a transformation law with
respect to an automorphy factor. The passage between the two descriptions
is just a choice of trivialization of the Hodge line bundle, and our choice corresponds to $\omega(z).$

The line bundle $\mathcal L^{\otimes k}$ has its natural Hermitian
metric induced by the bilinear form on the period line, the Petersson norm. In Yoshikawa's
notation this is written, after the same trivialization, as
\begin{equation}\label{eq: petersson norm and L2 norm}
        \|F\|^2_{\operatorname{Pet}}(z)
        =|F(z)|^2 (\omega(z), \overline{\omega(z)})^k = |F(z)|^2   (2(y,y))^k.
\end{equation}
For $k=1$ this is, up to factors, the usual $L^2$-norm of a section of $\mathcal{F}^2$.

\subsection{Application to Enriques surfaces}\label{subsect: Enriques dist sections} We specialize the above discussion to the setting of Enriques surfaces. We let $X$ be an K3 surface with a fixed point-free non-symplectic involution $\iota$ and $p:X \to Y=X/\iota$ the quotient map to the corresponding Enriques surface. 

\subsubsection{} We consider the Enriques lattice, $M=U(2) \oplus E_8(-2) \subseteq \Lambda_{K3}$, with  map given by  $(u,e_8) \mapsto (u, u, 0, e_8, e_8)$. The quotient is $U\oplus U \oplus E_8(-1)$ and hence $M$ is primitive. It is isomorphic to the image $p^\ast: H^2(Y, \Z)/\hbox{torsion} \subseteq H^2(X, \Z) \simeq \Lambda_{K3},$ which moreover identifies $M$ with the invariant subspace $H^2(X, \Z)^+$. Then $M^\perp = U \oplus U(2) \oplus E_8(-2)= H^2(X, \Z)^{-}$. For a detailed discussion about these embeddings we refer to \cite[Subsect 1.3. and 1.5.]{NamikawaEnriques}. 

There are two 0-cusps that arise by considering $M^\perp=U\oplus K_1$ or $M^\perp =U(2)\oplus K_2$, the level 1 and level 2 cusps. For the purposes of this article we focus on the level 1 cusp. The lattice $K_1$ also does not have any roots, since $\delta^2$ is always divisible by 4, so we can pick $C_K^+$ as our Weyl chamber. 

Any point of $\mathcal{D}^\circ_{M^\perp}$ defines a K3 surface that carries a non-symplectic involution which,  by results of Nikulin (\cite[Theorem 4.2.2]{NikulinFactorGroups}), is fixed point free. We fix a small  cone $\mathfrak{s} \subseteq M_\R$. Then Proposition \ref{prop:HodgeTatedegeneration} hence gives us a family of Enriques surfaces over $(\D^\ast)^{10}$, whose covering K3 surface has a strongly unipotent degeneration. 

\subsubsection{} On the Enriques moduli space one can define a special orthogonal modular form, the so called Borcherds product $\Phi$ \cite{Borcherdsgrassmannian, BorcherdsEnriques}. From equation \ref{eq:yoshikawaBCOV} of the introduction\ref{section:latticepolarizedk3},  Yoshikawa's results in \cite{YoshikawaK3inventiones}  give an expression of $\tau_{BCOV}$ . Combining this result with \eqref{eq: petersson norm and L2 norm}, this gives that the $|DRR|$ function in Definition \ref{def:DRR function} for Enriques surfaces can be expressed as:
\begin{equation}
    |DRR(q)|^2=\frac{\left(\tau_{BCOV}\cdot ||1||^2_{L^2,BCOV}\right)^{12} }{\left(||\omega(z)||^2_{L^\frac{2}{2}}\right)^{\frac{12}{2}}}=(C')^{12}|\Phi|^6
\end{equation}
where $C'$ is independent of $Y$ encoding both the constant $C$ from equation \ref{eq:yoshikawaBCOV} and the  $2\pi$-factor in the $L^2$-BCOV metric which is constant by \ref{table: DRR surfaces} . For comparison notice that this volume factor didn't appear in Yoshikawa's work.
  
  This result is of an equality between functions defined on the moduli space of Enriques. In order to talk about mirror symmetry we need to choose a cusp with respect to which we expand the Borcherd's product. 
  
\subsubsection{} Let us firstly say something about the different names for cusps appearing in the literature, on the preceding  parts of the paper we followed the notation of \cite{YoshikawaK3inventiones,YoshikawaTrinity}. In \cite{MR2391190} Klemm-Mariño argue that, from the point of view of physics, the  level $2$ cusp in \cite[Example 13.7]{Borcherdsgrassmannian} is the appropriate cusp where enumerative mirror symmetry is expected. This so called level $2$ for Borcherds corresponds to level $1$ for us. This is because in the above mentioned example Borcherds expands $\Phi$ with respect to the dual lattice of our $M^\perp$ and this has the effect of interchanging the cusps. With this in mind Borcherds obtains: 
$$\Phi(q)=\prod_{\lambda\in U\oplus E_8 }\left(\frac{1-q^\lambda}{1+q^\lambda}\right)^{c\left(\frac{\lambda^2}{2}\right)}$$
with $c(n)$ being the $n$-th coefficient of the Taylor expansion of a certain modular form he denotes as $f_{00}+f_{11}$, which could have an expansion with fractional exponents. In our special context the lattice over which we sum is even, so only coefficients $c(n)$ with $n\in \mathbb{N}$ matter. This means that we look at the expansion of $f_{00}=8\frac{\eta(2\tau)^8}{\eta(\tau)^{16}}=8\cdot \prod_n\left(\frac{1+q^n}{1-q^n}\right)^8$. 

\subsubsection{}On the other hand, in recent work \cite{Oberdieck2025KMformula} Oberdieck proves the so called Klemm-Mariño formula conjectured in \cite{MR2391190} by looking at the analytic torsion of $K3$s. Oberdieck computes the Gromov Witten theory of an Enriques surface and also obtains the genus-one generating function in terms of a Borcherds product. To be more precise Oberdieck's \cite[Corollary 1.2.]{Oberdieck2025KMformula} states that 
$$\exp\left(\sum N_{1,\beta}q^\beta\right)=\prod_{\beta\in H_2,\;\beta>0}\left(\frac{1+q^\beta}{1-q^{\beta}}\right)^{a\left(\frac{\beta^2}{2}\right)}$$
 where $a(n)$ the $n$-th Taylor coefficient of $\prod_n\left(\frac{1+q^n}{1-q^n}\right)^8$. This means exactly that 
 $$\exp(\sum N_{1,\beta}q^\beta)=\Phi^{\frac{-1}{8}}$$
When the Borcherds product is expanded near the level $1$ cusp. Putting everything together one gets that $-F_{1}^B(q)=F_{1}^A$ as desired.

  \subsubsection{} If one expands the above functions around the other cusp there is no apparent relation between $F_1^A$ and $F_1^B$. Geometrically, we chose to interpret it as the level one cusp avoiding an orbifold limit. This is because in \cite{alexeev2023compactificationsmodulispacesk3} it is shown that around the level $2$-cusp one can choose a one-dimensional Kulikov degeneration such that the involution extends with fixed points. On the other hand, around the level one cusp in \textit{loc.cit.} a one-dimensional family where the involution extends without fixed points is constructed.
  
Recall that our choice of a regular cone \(\sigma\) gives a toric partial compactification \(T_{K,\sigma}\) of \(T_K\), containing a polydisk \(\D^{10}\) with \((\D^\ast)^{10}\subset T_K\). By the Borel extension theorem, the period map
\[
(\D^\ast)^{10}\longrightarrow \Mcal_{\mathrm{Enr}}
\]
to the moduli space of Enriques surfaces extends to a holomorphic map
\[
\D^{10}\longrightarrow \Mcal_{\mathrm{Enr}}^{BB}
\]
of its Baily--Borel compactification.

The coordinate hyperplanes are mapped to the level-one cusp by construction. Indeed, the corresponding monodromy logarithms \(N_{k_i}\) satisfy \(N_{k_i}^2\neq 0\), and their monodromy weight filtrations have \(W_0=\mathbb Qe\), precisely the rational isotropic line defining the level-one cusp.

Consequently, because of the result of \cite{alexeev2023compactificationsmodulispacesk3}, any moduli interpretation of the boundary of our polydisk that is compatible with the Baily--Borel compactification must associate these boundary components with the level-one cusp. This motivates part of our discussion of cusps in the introduction (see \textsection \ref{subsubsect:cuspproblem intro}) and guides our subsequent choices of cusps in the case of bielliptic surfaces.

\section{Mirror constructions for bielliptic surfaces}\label{sect: geom biell}

In this section, we outline the construction of the families that we will use to prove our mirror symmetry statements for bielliptic surfaces. As we recall further down, see \ref{tab:bielliptic-classification}, the bielliptic surfaces come in seven different types, according to the group $G=\Z/N\Z 
\times \Z/K\Z$ which acts on a product of elliptic curves.  We will prove the following theorem:

\begin{theorem}\label{Thm: mirror bielliptic}
    For each bielliptic type, there is a bielliptic mirror family 
    $X \to \D$,
   given by \[X_q=(E\times F_q)/G,\;q\neq 0\] for a fixed elliptic curve $E$ and $F_q=(\C^\ast/q^\Z)$. More precisely, it satisfies the following properties:
    \begin{enumerate}
        \item\label{item1: mirror biell} The monodromy of $X_q$ is strongly unipotent in the sense of \ref{subsubsect: strongly unipotent}. 
        \item\label{item2: mirror biell} the singular fiber $X_0$ is constructed through a maximal N\'eron polygon for the family $F_q$, and $G$ acts without fixed points.
        \item\label{item3: mirror biell} The mirror map for this family is $q\mapsto q'=q^{\frac{1}{N}}$.
        \item\label{item4: mirror biell} The distinguished section in the sense of \ref{def:Hodge Tate distinguished sections}  is 
        $$\eta_{BCOV}=\left(1_{2,2}\otimes \frac{(2\pi)^2 i}{2\im(\tau')}d\overline{z}_{\tilde{F}}\right)^2.$$
        Here $e^{2\pi i \tau'}=q'$ and we have identified $H^1(X_q)=H^1(\tilde{F}_q)$ for $\tilde{F}_q=F_q/G$.
    \end{enumerate}
\end{theorem}
\begin{remark}
    In the above the parameter $q$ is what later on \ref{subsubsect: construction mirror family biell} will be referred to as $q_0$ to highlight that corresponds to the $0$-cusp. 
\end{remark}
The theorem is the combination of the following subsections. We will recall the relationship between moduli of bielliptic surfaces and modular curves. Part of the issue will then be to identify the correct cusp which is strongly unipotent in the correct sense, and admits good coordinates around it.

\subsection{Geometry of bielliptic surfaces}

We will begin by recalling some standard facts about bielliptic surfaces and their moduli to then recast them in our situation. We base our presentation on  \cite{BagneraDeFranchis1908, Compactcomplexsurfaces,SerranoBielliptics} and refer to this for additional details on the general material.
\subsubsection{}

Bielliptic surfaces are smooth projective surfaces with $h^{0,2}=0, h^{1,0}=1$ and $K_X$ torsion of order $N > 1$. Equivalently, it is a minimal surface whose Albanese morphism is an elliptic fibration without singular fibers. Their Hodge diamond is as follows: 
        \[
    \begin{array}{ccccc}
      && 1 && \\
      & 1 && 1 & \\
      0 && 2 && 0 \\
      & 1 && 1 & \\
      && 1 &&
    \end{array}
    \]
All bielliptic surfaces are in fact quotients
\[
    X=(E\times F)/G,
\]
where \(E\) and \(F\) are elliptic curves and \(G\) acts freely. We use the
convention that the action on \(F\) is by translations, while the action on
\(E\) is given by a finite-order automorphism. Since $h^{1,0}=1$ and  the fibers of the natural map 
\begin{equation}\label{def:albanesemap} X = (E \times F)/G \to \tilde{F}=F/G
\end{equation} are connected, it follows that $F/G$ is the Albanese of $X$.

If
\[
    G \simeq \Z/N\Z \times \Z/K\Z,
\]
we write \(G=\langle g_N,g_K\rangle\), where \(g_N\) has order \(N\) and \(g_K\)
has order \(K\). The generator of the first factor, \(g_N\), acts on \(E\) by an automorphism of
order \(N\), and on \(F\) by translation by an \(N\)-torsion point. 

\subsubsection{} There are seven types of bielliptic surfaces. They are listed in the table below, as in \cite{BagneraDeFranchis1908}. Here
\[
    E_i=\C/(\Z+i\Z), \qquad
    E_\rho=\C/(\Z+\rho\Z), \qquad
    \rho=e^{2\pi i/3}.
\]
Moreover,  \(t_\varepsilon\) denotes translation by \(\varepsilon\) on \(E\), and
the corresponding generators act on \(F\) by translations \(t_a\) and \(t_b\).
More precisely, in the cyclic cases the generator acts on \(F\) by \(t_a\).
In the product cases, the two generators act on \(F\) by \(t_a\) and \(t_b\),
respectively. We use
\[
    \varepsilon_i=\frac{1+i}{2},\qquad
    \varepsilon_\rho=\frac{1-\rho}{3}.
\]
A point of \(F[n]\) is called primitive if it has exact order \(n\).

\begin{table}[htbp]
\centering
\renewcommand{\arraystretch}{1.25}
\begin{tabularx}{\textwidth}{@{}c c c c X@{}}
\toprule
Type & \(G\) & \(E\) & \(E\)-action & translation data on \(F\) \\
\midrule

\(a_1\)
& \(\Z/2\Z\)
& arbitrary
& \(g_N:z\mapsto -z\)
& \(g_N=t_a \text{ for }a\in F[2]\setminus\{0\}\)
\\

\(a_2\)
& \(\Z/2\Z\times \Z/2\Z\)
& arbitrary
& \(g_N:z\mapsto -z,\quad g_K=t_\varepsilon\)
& \(g_N=t_a,g_K=t_b \text{ for }a,b\in F[2]\) independent
\\
\(b_1\)
& \(\Z/3\Z\)
& \(E_\rho\)
& \(g_N:z\mapsto \rho z\)
& \(g_N=t_a \text{ for }a\in F[3]\) primitive
\\

\(b_2\)
& \(\Z/3\Z\times \Z/3\Z\)
& \(E_\rho\)
& \(g_N:z\mapsto \rho z,\quad g_K=t_{\varepsilon_\rho}\)
& \(g_N=t_a,g_K=t_b \text{ for }a,b\in F[3]\) independent
\\
\(c_1\)
& \(\Z/4\Z\)
& \(E_i\)
& \(g_N:z\mapsto iz\)
& \(g_N=t_a \text{ for }a\in F[4]\) primitive
\\

\(c_2\)
& \(\Z/4\Z\times \Z/2\Z\)
& \(E_i\)
& \(g_N:z\mapsto iz,\quad g_K=t_{\varepsilon_i}\)
& \(g_N=t_a,g_K=t_b \text{ for }a\in F[4]\) primitive and  \(b\in F[2] \setminus \langle 2a\rangle\)
\\

\(d\)
& \(\Z/6\Z\)
& \(E_\rho\)
& \(g_N:z\mapsto -\rho z\)
& \(g_N= t_a\) for \(a\in F[6]\) primitive
\\

\bottomrule
\end{tabularx}
\caption{Classification table of the seven types of bielliptic surfaces.}
\label{tab:bielliptic-classification}
\end{table}
\FloatBarrier
A different presentation of these types can be given by considering the variety $A$ from Lemma \ref{lemma: covering associated to torsion KX} associated to the equality $K_X^{\otimes N}=\mathcal{O}_X$, by the classification of compact complex surfaces, is an abelian surface. In this case $X=A/(\Z/N\Z)$, we will use this alternative presentation to study the coarse moduli of bielliptic surfaces.
$A$ is a product for the cases $a_1,b_1,c_1,d$ so one refers to them as split bielliptic surfaces. For more details, such as an explicit presentation of $A$ or of the action of $\Z/N\Z$, see \cite{SuwaHyper,ModuliHyperell}.

\subsubsection{}
We want to study the integral variation of Hodge structures associated to a family of bielliptic surfaces. In general, by the K\"unneth formula $H^2(X)$ is given by $G$-invariant elements of $H^2(E)$  and $H^2(F)$. Since these are trivial VHS all the non constant part is in $H^1(X)$. 
 Moreover, the Albanese map in \eqref{def:albanesemap} induces an isomorphism on the level of  VHS associated to $H^1$. Hence all VHS is described by the cohomology of the Albanese $H^1(\widetilde{F},\Z)$. 

In contrast to past examples of genus 1 mirror symmetry(\cite{eriksson2022genus,Fang-Lu-Yoshikawatorsionthreefolds}) this is, up to our knowledge, the first example in genus one mirror symmetry where the information of the $B$-model is not encoded entirely in the middle cohomology group.

\subsection{Construction of the mirror family} \label{subsect: Moduli of biell}

It is shown in the works of Suwa \cite{SuwaHyper} and Tsuchihashi
\cite{ModuliHyperell} that the coarse moduli space of complex bielliptic
surfaces of a fixed type is described by a modular curve, or a product of two them depending on the type.

As $H^1(X,\Z)=H^1(\tilde{F},\Z)$ is an equality of VHS, the degeneration only depends on the parameter controlling the Albanese. We will see in \ref{vanishing for n,m not zero} that $\tau_{BCOV}$ only depends on the same one--parameter. Thus, we actively ignore the possible complex deformations from $E$ in the $a_1,a_2$ cases, and focus on a one dimensional moduli space.

\subsubsection{} The various \(N\)-torsion translations, and in particular the corresponding $N$-torsion points, of the elliptic curve \(F\), have as coarse moduli the level $N$ modular curve
\[
    Y_1(N)=\Gamma_1(N)\backslash \mathbb H.
\]
Here
\[
    \Gamma_1(N)
    =
    \left\{
    \begin{pmatrix}
    a & b\\
    c & d
    \end{pmatrix}
    \in \SL_2(\mathbb Z)
    \ \middle|\
    a\equiv d\equiv 1 \!\!\pmod N,\quad
    c\equiv 0 \!\!\pmod N
    \right\}
\]
with compactification $X_1(N).$

In order to construct the canonical coordinates, we need a description in analytic terms. If
\[
    F_\tau=\mathbb C/(\mathbb Z+\tau\mathbb Z),
    \qquad \tau\in\mathbb H,
\]
then we may choose the marked torsion point to be
\[
    P=\frac{1}{N}\in F_\tau[N].
\]
Thus the universal analytic picture is given by the family of elliptic
curves \(F_\tau\) together with the distinguished \(N\)-torsion point
\(1/N\), modulo the action of \(\Gamma_1(N)\).

\subsubsection{} 
The various cusps of $Y_1(N)$ all give natural candidates for mirror families of our bielliptic surfaces. The degenerations of the actual elliptic curves are presented as
$F_\tau=F_q=\C^\ast/q^\Z $ with $q = \exp(2\pi i\tau)$, degenerating into a $k$-gon. We make the choice of the 0-cusp of width $N$, and argue below, in particular in Proposition \ref{prop:Ngon}, why this is a natural choice which also satisfies the conclusions of Theorem  \ref{Thm: mirror bielliptic}. 

Since we are interested in the behavior of the Hodge structures  $H^1(X, \Z)$, it is natural to choose to degenerate simultaneously $X$ and $\tilde{F}$, relative to the Albanese map  \eqref{def:albanesemap}. We degenerate $F_q$ as $q \to 0$, this will automatically be strongly unipotent, which is a first step towards the first property of the mirror family construction ( \ref{Thm: mirror bielliptic}). The  degeneration of $X$ we are using is that of \cite[p.11]{ModuliHyperell} and it is related to that of $\widetilde{F}$ via the Albanese map
\begin{equation}\label{eq:albanesemorphism}
X_q=(E\times F_q)/G \to \tilde{F}=F_q/G \end{equation}
It seems natural, and in fact turns out to be necessary, to choose the degeneration such that the group action also behaves well on the singular fiber. This means avoiding an orbifold limit, so we degenerate $E \times F_q \to F_q$ in such a way that $G$ continues to act fixed point freely also on the singular fiber. 

We now explain in  \ref{prop:Ngon} how our choice of degeneration is equivalent to choosing a cusp in $X_1(N)$ given by a N\'eron $N$-gon (\cite{DeligneRapoport}), which corresponds, in the standard terminology, to the $0$-cusp of width $N$.

\begin{proposition}\label{prop:Ngon}
    The group action on $E \times F_q \to F$ extends to the singular fiber fixed point freely only around the $0$-cusp associated to an $N$-gon.
\end{proposition}
\begin{proof}
    
    The canonical coordinates around the different $j$-gons are associated to base changing by $z\mapsto z^j$ (see the following paragraphs where we explicitly use this to construct the family or \cite{DeligneRapoport,DiamondShurman}). This means that over the $N$-gon the group action $\tau\mapsto \tau+\frac1N$ consists of permuting irreducible components associated to different branches of the $N$-th root. For other cusps the action has fixed points. 
\end{proof}

\subsubsection{}\label{subsubsect: construction mirror family biell} We now recall the construction of a family over $\H$ given by $A_\omega/\left(\Z/N\Z\right)$ where  $A_\omega$ is the abelian variety $A_\omega=\C^2/\Z^4\Omega$ with periods 
\begin{equation}\label{eq:periods abelian cover}
    \Omega=\begin{pmatrix}
        \alpha&0 \\
        p_0& \omega\\
        1&0\\
        0&1
        
    \end{pmatrix}
\end{equation}

Here $p_0,\alpha$ are determined by the type of bielliptic surface, in the notation of table \ref{tab:bielliptic-classification} $p_0=\varepsilon$ and $\alpha$ is the parameter of $E$. 
This corresponds to a regular fiber $X_\tau=(E\times F_\tau)/(\Z/N\Z\times \Z/K\Z)$. Notice that we have a degree $K$ cover $E\times F_\tau\to A_\omega$  which gives $\omega=\frac{\tau}{K}$ (we call $\omega$ what \cite{ModuliHyperell} calls $\tau$).

To study the $0$ cusp of width $N$, we perform the change of variables $\omega'=\frac{-1}{\omega}$, exchanging the 0 and $\infty$ cusps. This induces the isomorphism:
\begin{equation*}
    \begin{split}
        \left(\C/(\Z+\omega\Z),\frac{1}{N}\right)&\to\left(\C/(\Z+\omega'\Z),\frac{\omega'}{N}\right), \quad z \mapsto \frac{-1}{\omega} \cdot z.
    \end{split}
\end{equation*}

The natural parameter around this point of the moduli space is  $q'=q_{0}^\frac{1}{N}=e^{2\pi i \frac{-1}{\omega}\frac{1}{N}}=e^{2\pi i \frac{\omega'}{N}}$, see \cite[page 60 ]{DiamondShurman}, notice that $q_{0}$ is the parameter after exchanging the cusps i.e. $q_{0}=e^{2\pi i \omega '}$. Then  $q'=q_{0}^\frac{1}{N}=e^{2\pi i \tau'}$ gives naturally a coordinate over the disk centered at the point we degenerate to. In terms of $\tau$ this corresponds to $\tau'=\frac{-K}{\tau N}$.

\subsubsection{}

Notice that, with respect to this parameter, the geometric monodromy 

around $q'=0$ corresponds to $\tau' \mapsto \tau' +1$ on $\tilde{F}$. The action on the homology cycles can be described by picking the basis $\gamma_0,\gamma_1$ on $H_1(\tilde{F},\Z )$ where $\gamma_0$ corresponds to the segment connecting $0$ and $1$ and $\gamma_1$ connects $0$ and $\tau'=\frac{-K}{N\tau}$. Then we have $T(\gamma_0)=\gamma_0$ and $T(\gamma_1)=\gamma_0+\gamma_1$ where $T$ is the monodromy operator in homology (see \cite[Eq 1.12 p.18]{P.Mapings_P.domains_2017}).

Hence, following  Deligne \cite{DeligneLocalBehav} and Morrison \cite[Section 3]{Morrisonguide}, we have that the mirror map is a quotient of the periods of $\gamma_1$ by that of $\gamma_0$. We hence propose the following definition:\\
\begin{definition}\label{def: Mirror map bielliptic}
    The mirror map for our mirror family $(F_\tau\times E)/G$ around the $0$-cusp  is given by $\tau \mapsto \tau'=\frac{-K}{\tau N}$. Equivalently it is given by $\tau \mapsto q_{0}^\frac{1}{N}$
\end{definition}

\subsubsection{}\label{subsubsect: (bi)elliptic distinguished section}
We now finalize with the computation of the distinguished section $\eta_{BCOV}$. In terms of homology we have $W_{-2}H_1=\langle \gamma_0\rangle\cong\left(H^1/W_0H^1\right)^\vee$ and  $W_0/W_{-2}$ is generated by $\gamma_1$. By \ref{def:Hodge Tate distinguished sections}, a distinguished section $\eta_0$ in  $\mathcal{F}^0/\mathcal{F}^1$ can uniquely be written as a lift to a section $\eta \in \mathcal{F}^0$ such that:
\begin{displaymath}
    \int_{\gamma_0}\eta=0, \quad \int_{\gamma_1}\eta=(2\pi)^2.
\end{displaymath}

We write $\eta= a dz + b d\overline{z}$ and compute $a,b$ from
$a\int_{\gamma_0}dz+b\int_{\gamma_0}d\overline{z}=0$ and $a\int_{\gamma_1}dz+b\int_{\gamma_1}d\overline{z}=(2\pi)^2$.
We find $b=-a$ and $a=\frac{-(2\pi)^2i}{2\im(\tau')},$
and projecting this $\eta$ onto $\mathcal{F}^0/\mathcal{F}^1\cong H^{0,1}$ we find that $$\eta_0=\frac{(2\pi)^2 i}{2\im(\tau')}d\overline{z}_{\tilde{F}}$$ 
We conclude that 
$$\eta_{BCOV}=(1_{2,2}\otimes \eta_0)^2,$$
which also finishes the proof of Theorem \ref{Thm: mirror bielliptic} and concludes this section.

\section{$F_1^B$ of a bielliptic surface}\label{Sect: F_1B bielliptic}

The end goal of this section is to finalize the study of the $B$-side for bielliptic surfaces.  This amounts to the computation of $F_1^B$ from the mirror family of bielliptic surfaces obtained in Theorem \ref{Thm: mirror bielliptic}. Our main result is:

\begin{theorem}\label{Thm:F_1 bielliptic}
    Let $X\to \mathbb{D}^\ast$ be the constructed bielliptic mirror family with mirror coordinate $q'$. Then 

    $$F_1^B(q')=\sum_{j|N}\lambda(j)\log(\eta(q'^j))$$
    
     where $(\lambda(j))_{j|N}$ are given in the following table, the components of are ordered by increasing $j$.
     \[
     \renewcommand{\arraystretch}{1.4} 
\setlength{\tabcolsep}{10pt}
\begin{array}{|c|c|c|c|}
\hline
\text{Type } a\; (N=2) & \text{Type } b\; (N=3) & \text{Type } c \; (N=4)& \text{Type } d \; (N=6)\\
\hline
\begin{array}{c}
\lambda=(4,-2)
\end{array}
&
\begin{array}{c}
\lambda=(3,-1)
\end{array}
&
\begin{array}{c}
\lambda=(2,1,-1)
\end{array}
&
\begin{array}{c}
\lambda=(1,1,1,-1)
\end{array}
\\
\hline
\end{array}
\]

\end{theorem}
The general definition of $F_1^B$ is stated in \ref{def:DRR function}, this function appears as a logarithm of the function $|DRR|$ which is given as the quotient of the norms of two distinguished sections. The full computation of $|DRR|$ is in terms of the $BCOV$ invariant which we compute in Theorem \ref{thm:computation of BCOV invariant} by studying equivariant data on $E\times F $. 
\subsection{Description of the $L^2$ inner products}\label{subsect:L2 products}
\subsubsection{} We want to study the relation between $L^2$ norms in the first cohomology groups of the various manifolds appearing. In particular throughout the text we use Arakelov theoretic renormalizations, that is, locally the Kähler form is: 
$$\omega=\frac{i}{2\pi}\sum_{j,k}h_{j,k}dz_j\wedge d\overline{z}_k$$

 This means that the $L^2$-inner product is $(\alpha,\beta)_{L^2}=\int_X h(\alpha,\beta)\frac{\omega^n}{n!}=\frac{1}{(2\pi)^n}\int_X\alpha\wedge \overline{\ast}\beta$. In particular, we  normalize to volume $\frac{1}{(2\pi)^n}$. These renormalizations will be important in specific parts of the following section, nonetheless the reader can assume we are in a volume $1$ situation except when we talk about metric considerations related to the Grothendieck-Riemann-Roch isomorphism. We identify $\mathcal{A}^{p,q}(X)=(\mathcal{A}^{p,q}(E\times F))^G$ as vector spaces. The Kähler form $\omega_{E\times F}$ is $G$-invariant so it induces a form $\alpha\in \mathcal{A}^{1,1}(X)$ such that $\pi^\ast\alpha=\omega_{E\times F}$.

We have $\int_{E\times F}\omega_{E\times F}^2=\int_{E\times F}\pi^\ast\alpha^2=|G|\int_X\alpha^2$. Thus, the Kähler form on $X$ according to our normalizations is $\omega_X=\sqrt{|G|}\alpha$. 

\subsubsection{} In order to compute the factor $B(X)$ we first establish the following Lemma comparing the relevant $L^2$ norm on the bielliptic surface to that of the elliptic curve $\tilde{F}$ (recall the morphism \ref{def:albanesemap}):

\begin{lemma}\label{relation between L2 norm of surface and curve}
    Let $p:X\to \tilde{F}$ be the Albanese morphism.    
    The following equality holds for the $L^2$ norms on $H^1(X)$ and $H^1(\tilde{F})$ chosen with our previously stated normalizations:
    $$(p^\ast(\alpha),p^\ast(\beta))_{L^2,X}=\frac{1}{2\pi}\sqrt{|G|}(\alpha,\beta)_{L^2,\tilde{F}}$$
\end{lemma}
\begin{proof}
Since $H^1(X)$ coincides with its primitive part we can actually describe the $L^2$ inner product as $$(p^\ast\alpha,p^\ast\beta)_{L^2}=-\frac{1}{2\pi}\int_{X}p^\ast\alpha\wedge C \overline{p^\ast\beta} \wedge \omega_X$$ 
here $C$ is the Weil operator. The above coincides with 
\begin{equation*}
    \begin{split}
        -\frac{1}{|G|2\pi}\int_{E\times F}\pi^\ast p^\ast\alpha\wedge C\overline{\pi^\ast p^\ast\beta}\wedge\sqrt{|G|}\pi^\ast(\alpha)&=\\-\frac{1}{|G|2\pi}\int_{E\times F}\pi^\ast p^\ast\alpha\wedge C\overline{\pi^\ast p^\ast\beta}\wedge\sqrt{|G|} \omega_{E\times F}&=\frac{1}{\sqrt{|G|}}(\pi^\ast p^\ast \alpha,\pi^\ast p^\ast \beta)_{L^2,E\times F}.
    \end{split}
\end{equation*}
Here $\pi: E\times F\to X$ and $\phi: F\to \tilde{F}$ are the natural quotient maps, which are both degree $|G|$ covers.

We now furthermore have the following
\begin{equation*}
    \begin{split}
        (\pi^\ast p^\ast \alpha,\pi^\ast p^\ast \beta)_{L^2,E\times F}=-\frac{1}{2\pi}\int_{E\times F}\pi^\ast p^\ast\alpha\wedge C\overline{\pi^\ast p^\ast\beta}\wedge(p_1^\ast\omega_E+p_2^\ast \omega_F)&=\\       
        \left( \int_E\omega_E\right)\left(-\frac{1}{2\pi}\int_F\phi^\ast(\alpha\wedge C\overline{\beta})\right)=-\frac{1}{2\pi}|G|\frac{1}{2\pi}\int_{\tilde{F}}\alpha\wedge C\overline{\beta}&=\frac{|G|}{2\pi}(\alpha,\beta)_{L^2,\tilde{F}}.
        \end{split}
\end{equation*}
From this the result follows.
\end{proof}
 The map $p^\ast: H^1(\widetilde{F}, \Z) \to H^1(X, \Z)$ is an isomorphism compatible with families. Remembering the factor relating both $L^2$-norms we will often identify  these two polarized variations of Hodge structures.

\subsubsection{} With this comparison at hand we can compute the first piece of the BCOV invariant; the factor $B(X)$. 
\begin{lemma}\label{lemma:computation B(biell)}
    Let $X=(E\times F_\tau )/G$ be a bielliptic surface and let $B(X)$ be as in \ref{def:B-term}. Then 
    $$B(X)=\frac{1}{|G|}$$
\end{lemma}
\begin{proof}
     By item \ref{item3:renormalized L^2 norm} in  \ref{prop:renormalized L^2 norm} we compute $$B(X)=\vol_{L^2}(H^0(X,\Z))^2\vol_{L^2}(H^1(X,\Z))^{-1}$$ 
    
    This first factor is equal to $\frac{1}{(2\pi)^4 }$, since $\vol(H^0)=\operatorname{covol}(H^0)^2=\vol(X)$ , as we renormalized our hermitian product at the beginning to be of volume $\frac{1}{(2\pi)^2 }$.

    Applying Lemma \ref{relation between L2 norm of surface and curve} we have $\vol_{L^2}(H^1(X,\Z))=\left(\frac{1}{(2\pi)}\sqrt{|G|}\right)^{b_1(X)}\cdot\vol_{L^2}(H^1(\tilde{F},\Z))$. Poincar\'e duality on $\tilde{F}$ gives $H^1(X,\Z)\times H^1(X,\Z)\to \Z(-1)$. This map is compatible with Serre duality, which is an $L^2$-isometry, so the volume of the middle cohomology is $\vol_{L^2}(H^1(\tilde{F}),\Z)=\frac{1}{(2\pi)^2}$. This proves the theorem.
\end{proof}

\subsubsection{}\label{L^2 part of DRR biell}We can now explicitly compute all of the $L^2$ contributions on the definition of $F_1^B$:  From Theorem \ref{Thm: mirror bielliptic} we saw that the canonical section is $$\eta_{BCOV}=(1_{2,2}\otimes \eta_0)^2,$$
where $\eta_0=\frac{(2\pi)^2 i}{2\im(\tau')}d\overline{z}_{\tilde{F}})^2$. It has $L^2$ norm $||\eta_{BCOV}||_{L^2}^2=||1||^4_{L^2}||\eta_0||^4_{L^2,X}=\frac{1}{(2\pi)^4}\frac{|G|}{(2\pi)^2}(\eta_0,\eta_0)^2_{L^2,\tilde{F}}.$
 The latter factor can be made explicit:
\begin{equation*}
    \begin{split}
        (\eta_0,\eta_0)_{L^2,\tilde{F}}&=-\frac{1}{(2\pi)}\int_{\tilde{F}}\eta_0\wedge C\overline{\eta_0}
        =\frac{1}{2\pi}\frac{-(2\pi)^4}{(2\im(\tau'))^2}\int_{\tilde{F}}-idz\wedge d\overline{z}\\
        &=\frac{(2\pi)^3}{2\im(\tau')}
    \end{split}
\end{equation*}

Thus,\begin{equation}\label{eq: quillen metric dist section bielliptic}
    \begin{split}
        ||\eta_{BCOV}||^2_{Q,BCOV}&=A(X)\cdot \mathcal{T}_{BCOV}\cdot (\eta_{BCOV},\eta_{BCOV})_{L^2}=|G|\frac{1}{4\im(\tau')^2}\mathcal{T}_{BCOV}
    \end{split}
\end{equation}
Since the metric is Ricci flat and the Euler characteristic is $0$ the Proposition \ref{prop:QuillenBCOV} gives the second equality. This reduces the computation to that of $\mathcal{T}_{BCOV}$ and $\tau_{BCOV}$. This is accomplished in Theorem \ref{thm:computation of BCOV invariant} below. 

\subsection{The zeta function of a bielliptic surface}
The goal now is to compute holomorphic analytic torsion of $X$ in terms of $G$-invariant forms on $E\times F$. On the product, the Laplacian and its eigenforms are well known and admit an explicit description which we now briefly recall.

\subsubsection{} Let us denote by $\Delta$ (without upper indices) the usual Laplacian on the plane times $-1$, that is, $\Delta=-(\frac{\partial^2}{\partial x^2}+\frac{\partial^2}{\partial y^2})$. On a single elliptic curve $C=\mathbb{C}/(\Z+\tau \Z)$ we take the volume form $\frac{i}{(2\pi)2\im(\tau)}(dz\wedge d\overline{z})=\frac{1}{2\pi \im(\tau)}(dx\wedge dy)$. This renormalizes the different $(p,q)$ Laplacians to be:
\begin{equation}\label{laplacians}
    \Delta^{p,q}_{\overline{\partial}}(f\cdot \alpha)= \frac{\im (\tau)}{2}\Delta(f)\cdot \alpha
\end{equation}
where $\alpha$ is among $1,dz,d\overline{z},dz\wedge d\overline{z}$. We also define 
\begin{equation}\label{def:enm} e_{n,m}(x,y)=e^{2\pi i nx}e^{2\pi i\frac{m-n\tau_1}{\tau_2}y}
\end{equation}
It is a an immediate result that the set $\{e_{n,m}\}$ forms a basis of eigenfunctions for $\Delta$.

In general we want to compute $\det(\Delta^{p,q}_{\overline{\partial}})$ in order to then compute $\mathcal{T}_{BCOV}$. This is defined as  $$\det \Delta^{p,q}=\exp\left(-\frac{d}{ds}\zeta^{p,q}(s)|_{s=0}\right)$$
where $\zeta^{p,q}(s)=\sum_{\lambda \text{ eigenvalue }}\frac{1}{\lambda^s}.$ For completeness we record the following standard lemma.

\begin{lemma}\label{lemma: zeta function elliptic curve}
    The regularized zeta function on the elliptic curve $C=\C/\Z+\tau\Z $ associated to $\Delta_{\overline{\partial}}^{p,q}$ is
    \begin{equation}\label{zeta funtion as E.series}
        \zeta^{p,q}(s)=E(\tau,s)=\frac{1}{(2\pi^2)^{s}}\sum_{m,n}\frac{\im(\tau)^s}{|m+n\tau|^{2s}}
    \end{equation}
    Notice the factor of $\frac{1}{(2\pi^2)^s}$ in our definition of the Einsenstein series.
\end{lemma}

\subsubsection{}\label{subsubsec:Kronecker} We state and prove the following Lemma, which is standard. However, in these formulas there are notorious discrepancies in the literature with respect to the constants. 
\begin{lemma}\label{lemma:torsion ell.curve}
    Let $E=\C/(\Z+\tau\Z)$ be an elliptic curve and $T=T(E,\mathcal{O}_E)$ the holomorphic analytic torsion of its trivial line bundle. Then:
    \begin{enumerate}
        \item\label{item1:torsion ell.curve} $\mathcal{T}_{BCOV}=T^{-1}$
        \item\label{item2:torsion ell.curve} $T=\frac{1}{2}\frac{1}{\im(\tau)}|\Delta(\tau)|^{\frac{-1}{6}}$ where $\Delta(\tau)=q\prod_n(1-q^n)^{24}$;  $q=e^{2\pi i \tau}$.
    \end{enumerate}
\end{lemma}
\begin{proof}
    The proof of \ref{item1:torsion ell.curve} follows from the definitions which in this case amount to $$T=(\det\Delta^{0,1})^{-1}= (\det\Delta^{1,1})^{-1}$$ 
    where the equality of determinants follows from \ref{lemma: zeta function elliptic curve} 
    and $$\mathcal{T}_{BCOV}=\det\Delta^{1,1}.$$

    To prove statement \ref{item2:torsion ell.curve} we apply Kronecker's limit formula at $s=0$ for the series $\sum_{m,n}\frac{1}{|m+n\tau|^{2s}}$. In particular we use \cite[Chapter VIII equation 17]{EllipticacordingtoEinsensteinandKronecker} which gives $$\frac{d}{ds}\left(\sum_{m,n}\frac{1}{|m+n\tau|^{2s}}\right)|_{s=0}=-\frac{1}{12}\log|(2\pi)^{12}\Delta(\tau)|^2.$$ 
    Since $$T=\exp\left(\frac{d}{ds}\left(\left(\frac{\im(\tau)}{2\pi^2}\right)^s\sum_{m,n}\frac{1}{|m+n\tau|^{2s}}\right)_{s=0}\right)$$ the result follows. 
\end{proof}

 We notice in passing that the Kronecker limit formula from \cite{EllipticacordingtoEinsensteinandKronecker} gives 
 $$\sum\frac{1}{|m+n\tau|^{2s}}=-1-\frac{s}{12}\log( |(2\pi)^{12}\Delta|^2)+O(s^2)$$
 which already gives $E(\tau,0)=-1$, we will use this fact later.

\subsubsection{} In the following subsections we will compute the spectral zeta functions corresponding to the bielliptic surface, by comparing with the product of the elliptic curves. All the computations that we do are summed up in the following result: 

\begin{theorem}\label{thm: global expression of combination of zeta functions}
    Let $X=(E\times F)/(\Z/N\Z\times \Z/ K\Z)$ be a bielliptic surface and $\tau'$ its mirror coordinate as defined in \ref{def: Mirror map bielliptic}. The following equality holds:
    
   \[ \sum_{p,q}(-1)^{p+q}pq\zeta^{p,q}(s)=\frac{1}{K^s}\left(\sqrt{|G|}\right)^s\left(\sum_{j|N}\frac{1}{\left(\frac{N}{j}\right)^s}\lambda(j)E(j\tau',s)\right)\]
   Here $(\lambda(j))_{j|N}$ has the following forms where the entries are ordered in increasing values of $j$. 
   \[
\begin{array}{|c|c|c|c|}
\hline
\text{Type } a & \text{Type } b & \text{Type } c & \text{Type } d \\
\hline
\begin{array}{c}
(4,-2)
\end{array}
&
\begin{array}{c}
(3,-1)
\end{array}
&
\begin{array}{c}
(2,1,-1)
\end{array}
&
\begin{array}{c}
(1,1,1,-1)
\end{array}
\\
\hline
\end{array}
\]
\end{theorem}

The main qualitative observation about this combination of spectral zeta functions can be found in  Prop.\ref{vanishing for n,m not zero}, where we see that it is a function only on one moduli parameter for all types of bielliptic surfaces. This can a priori be surprising because for type $a$ the  complex moduli is two dimensional.

\subsubsection{}Recall the renormalizations made on the $L^2$-metrics and the relation between the metrics on $E\times F$ and $X$ as stated in \ref{subsect:L2 products}. This means that the Laplacian we are interested in on $E\times F$ is of the form $\frac{1}{\sqrt{|G|}}(\Delta_E+\Delta_F)$. 

Thus, its eigenvalues, whenever the eigenvectors exist, are $\frac{1}{\sqrt{|G|}}(\lambda(\tau_E;n,m)+\lambda(\tau_F,k,l))$. Here 
$$\lambda(\tau;n,m)=(2\pi^2)\frac{|m+n\tau|^2}{\im(\tau)}$$ denotes each term in the sum  \ref{zeta funtion as E.series}. Let us also introduce the following notation: 
\begin{itemize}
    \item $\zeta=\zeta^{2,2}=\zeta^{0,0}$. 
    \item $\zeta_1=\text{spectral zeta function of the eigenvalues of $(1,0)\otimes (0,0)$}$.
    \item $\zeta_2=\text{ spectral zeta function of the eigenvalues of $(0,1)\otimes (0,0)$}$.
\end{itemize}

We reduce the computation of the BCOV torsion to the following result.
\begin{lemma}\label{lemma: reduction of spectral combination}
    We have 
\begin{equation}\label{BCOV formula}
    \sum_{p,q}(-1)^{p+q}pq\zeta^{p,q}= 2\zeta-\zeta_1-\zeta_2.
\end{equation}
\end{lemma}
\begin{proof}
    The eigenforms on $E\times F$ naturally decompose according to eigenforms on $E$ and $F$ of bidegrees summing up to the original one. Hence the spectral zeta functions are naturally sums of zeta functions corresponding to bidegrees on $E$ and $F$. 
    
    From the fact that $dz_2,d\overline{z}_2$ are $G$-invariant we have that these refined spectral zeta functions built from these bidegrees are independent of the second bidegree. Moreover, $\omega_E$ is also invariant and so the zeta function of $(1,1)\otimes (0,0)$ and coincides with the one for $(0,0)\otimes (0,0)$ (i.e. with $\zeta^{0,0}$). 
    
    From this one easily concludes that  
    $\zeta^{1,1}=2\zeta+\zeta_1+\zeta_2, \zeta^{2,1}=\zeta+\zeta_1$ and $\zeta^{1,2}=\zeta+\zeta_2$, and the result follows from definition.
\end{proof}

\subsubsection{} To compute the eigenforms we use the following surjection 
$$\phi^{p,q}:\mathcal{A}^{p,q}(E\times F)\to (\mathcal{A}^{p,q})^G(E\times F)=\mathcal{A}^{p,q}(X)$$ 
$$\alpha\mapsto \frac{1}{|G|}\sum_{g\in G}g\cdot \alpha$$ 
By the above argument we just want to study this map for forms that are $(0,0),(1,0)\otimes (0,0),(0,1)\otimes (0,0)$. To see the contribution of $G$-invariant eigenvalues, it is enough to study this map on an eigenbasis of $E\times F$.

\begin{remark}\label{remark: fixing group action bielliptic }
    To do our explicit computations it is important to recall (see \ref{tab:bielliptic-classification}) that we choose to represent the bielliptic surface given by the following group action: let $g_N$ be the generator of the $N$-torsion part of $G$ and $g_K$ the generator of the $K$-torsion part. Then $g_N$ acts on $E\times F$ as $(z_1,z_2)\mapsto (e^{\frac{2\pi i}{N}}z_1,z_2+\frac1N)$ and $g_K$ acts on $E$ as transalation by $\frac{\tau_E+1}{2},\frac{1-e^{\frac{2\pi i}{3}}}{3}$ or $\frac{1+i}{2}$ in cases $a_2,b_2,c_2$  and on $F$ as $$z_2\mapsto z_2+\frac{\tau_F}{K}.$$
\end{remark}

We have the following relation between split and non-split cases. In this formulation we let $e_{n,m}$ be an eigenfunction on $E$ and $e_{k,l}$ be an eigenfunction on $F$, both defined as in \ref{def:enm}.
\begin{lemma}\label{lemma:phi p,q function on split and nonsplit}
    The projection $\phi^{p,q}$ to the $G$-invariant part is explicitly
    $$\phi^{p,q}(e_{n,m}e_{k,l}\alpha)=P_K\left(e^{\frac{2\pi i (m\pm n-l)}{K}}\right)\cdot \phi^{p,q}_{\text{split}}(e_{n,m}e_{k,l}\alpha)$$
    Here $\alpha\in \{1,dz_1,d\overline{z_1}\}$, $P_K$ is the polynomial 
    $$P_K(t)=\frac{1}{K}\left(1+t+t^2+\ldots +t^{K-1}\right)$$ 
    and $\phi_{\text{split}}^{p,q}$ is the function appearing in the respective split case. Here we write $m\pm n$ in the exponent inside $P_K$ specifically because it is $m+n$ in cases $a_2,c_2$ and $m-n$ in case $b_2$.
\end{lemma}
\begin{proof}
    Notice that if $X$ is already split, $K=1$, and $P_K=1$, so the lemma is only interesting in the non-split case. By definition the action of $\phi^{p,q}$ is given by the following element of the group algebra $$\frac{1}{|G|}\sum_{0\leq i\leq N,0\leq j\leq K}g_N^ig_K^j=\left(\frac{1}{N}\sum_{i=0}^Ng_N^i\right)\left(\frac{1}{K}\sum_{j=0}^Kg_K^j\right).$$ The first factor is the definition of $\phi_{\text{split}}^{p,q}$. We moreover have that $dz_1,dz_2,d\overline{z}_1,d\overline{z}_2$ are all invariant under $g_K$. Finally, by inspection of the table \ref{tab:bielliptic-classification} we can conclude that $g_K\cdot (e_{n,m}e_{k,l})=e^{\frac{2\pi i (m\pm n-  l)}{K}}\cdot e_{n,m}e_{k,l}$ we can conclude.
\end{proof}

\subsubsection{} As anticipated in \ref{subsect: Moduli of biell} we are now in position to prove that $\mathcal{T}_{BCOV}$ only depends on $\tau'$. 
\begin{proposition}\label{vanishing for n,m not zero}
    The contributions to equation \eqref{BCOV formula} from the eigenvalues with $(n,m)\neq (0,0)$ cancel out. In particular,  $\tau_{BCOV}$ only depends on the moduli parameter $\tau'$.

\end{proposition}
\begin{proof}
By Lemma \ref{lemma:phi p,q function on split and nonsplit} it suffices to check this for the split case.
    
When $(n,m)\neq 0$ fix an eigenvalue $\lambda=\frac{1}{\sqrt{|G|}}\left(\lambda(n,m)+\lambda(k,l)\right)$. By Lemma \ref{lemma: reduction of spectral combination} we only study bidegree $(0,0)\otimes(0,0)$, $(1,0)\otimes (0,0)$ or $(0,1)\otimes (0,0)$ We denote by $V_{0}$, $V_1$, $V_{-1}$ the respective $\lambda$-eigenspaces.
We have the following relations:
     \begin{itemize}
         \item $V_1\cong V_0\otimes \C_N$
         \item  $V_{-1}\cong V_0\otimes (\C_N)^{\otimes N-1} $
     \end{itemize} 
     where $\C_N$ is the representation where  $g_N$ acts as multiplication by $e^{\frac{2\pi i}{N}}$.  For a representation $(G,V)$ with character $\chi$ we have $\dim V^G = \frac{1}{|G|}\sum_{g \in G} \chi(g)$. It is straightforward to verify that $g_N^j(e_{n,m} e_{k,l})=C_{n,m}\cdot e_{n',m'} e_{k,l}$ with $(n,m)\neq(n',m')$ whenever $j \neq 0$. This means that the characters appearing vanish because the matrix representing $g_N^j$ has only zeroes along its diagonal. Thus we get the agreement for the sums of characters.
     Since the dimensions of these invariant spaces agree for the above representations, and  we deduce the cancellation of their corresponding contribution to \eqref{BCOV formula}. The second assertion follows from the fact that, up to constants, $\tau_{BCOV}$ is equal to the derivative at $0$ of the combination of zeta functions we are studying.

\end{proof}
\subsubsection{} We introduce the notation 
\[S(\tau,s)=\sum_{p,q}(-1)^{p+q}pq\zeta^{p,q}(s),\]
which allows us to write $\mathcal{T}_{BCOV}=\exp(-\frac{d}{ds}S(\tau,s)|_{s=0})$.
To finish our computation we study $(n,m)=(0,0)$. From Lemma \ref{lemma:phi p,q function on split and nonsplit} and a computation analogous to the one in its proof, we find the following lemma:
\begin{lemma}\label{lemma:phie_0,0}
The following equalities hold:
\begin{itemize} 
    \item $\phi^{0,0}(1\otimes e_{k,l})=P_K(e^{\frac{2\pi i l}{K}})P_N(e^{\frac{2\pi i k}{N}})1\otimes e_{k,l}$.
    \item $\phi^{1,0}(1\otimes e_{k,l}dz_1)=P_K(e^{\frac{2\pi i l}{K}})P_N(e^{ \frac{2\pi i (k+1)}{N}})1\otimes e_{k,l}dz_1$.
    \item $\phi^{0,1}(1\otimes e_{k,l}d\overline{z}_1)=P_K(e^{\frac{2\pi i l}{K}})P_N(e^{\frac{2\pi i (k-1)}{N}})1\otimes e_{k,l}d\overline{z}_1$.
\end{itemize}

The polynomials $P_N$ and $P_K$ happen to vanish at non-trivial roots of unity, and this only depends on the values of $l\text{ mod }K$ and $k \text{ mod } N$. Otherwise they are equal to 1. This means that $\phi^{p,q}$ either vanishes or is the identity.
\end{lemma}

The main idea of the computation of $S(\tau,s)$ is to express the sums of the type $\sum_{k=1 \text{ mod }N} \frac{1}{\lambda(k,l)^s}$ as $E(j\tau,s)$ for a certain $j$. We present with full detail only the cases $a_1,a_2$ and $b_1$. The rest are done in the same way. 

\begin{itemize}
    \item Case $a_1$: In this case $K=1$ and $N=2$. By Lemma \ref{lemma:phie_0,0} the expressions for which $\phi^{0,0} e_{k,l}$ is non-trivial correspond to even $k$. Similar considerations give that the other two don't vanish only on odd $k$. Hence we find that: $$S(\tau,s)=2\zeta-\zeta_1-\zeta_2=\left(\sqrt{|G|}\right)^s\left(2\sum_{k,l}\frac{1}{\lambda(2k,l)^s}-\sum_{k,l}\frac{1}{\lambda(2k+1,l)^s}-\sum_{k,l}\frac{1}{\lambda(2k-1,l)^s}\right)=$$
    $$\left(\sqrt{|G|}\right)^s\left( 4\sum_{k,l}\frac{1}{\lambda(2k,l)^s}-2\sum_{k,l}\frac{1}{\lambda(k,l)^s}\right).$$ 
    We now use that,  substituting the value of $\lambda(2k,l)$, the following holds:
    $$\sum_{k,l}\frac{1}{\lambda(2k,l)^s}=(2\pi^2)^{-s}\sum_{k,l}\frac{\im(\tau)^s}{|l+2k\tau|^{2s}}=    \frac{1}{(2\pi^2)^{s}2^s}\sum_{k,l}\frac{\im(2\tau)^s}{|l+k(2\tau)|^{2s}}=\frac{1}{2^s}E\left(2\tau,s\right).$$ 
   
    We may interpret this formula in terms of the mirror coordinate \ref{def: Mirror map bielliptic}  $\tau'=\frac{-1}{N\tau}$ to get: $$S(\tau',s)=\left(\sqrt{|G|}\right)^s\left(4\frac{1}{2^s}E\left(\tau',s\right)-2E\left(2\tau',s\right)\right)$$  
    Here we use the modular transformation property $E(\tau,s)=E(\frac{-1}{\tau},s)$ for any $\tau\in \H$.

    \item Case $a_2$: Firstly notice that we have the same parity conditions on $k$ as in the split case. Moreover, since we have the factor $P_2(e^{\frac{2\pi i l}{2}})$,  $l$ is necessarily even. So we get:
    $$S(\tau,s)=\left(\sqrt{|G|}\right)^s\left(2\sum_{k,l}\frac{1}{\lambda(\tau;2k,2l)}-2\sum_{k,l}\frac{1}{\lambda(2k+1,2l)^s}\right)=\left(\sqrt{|G|}\right)^s\left(4\sum_{k,l}\frac{1}{\lambda(2k,2l)^s}-2\sum_{k,l}\frac{1}{\lambda(k,2l)^s}\right). $$ 
    The first term in the parenthesis is equal to $\frac{1}{2^{2s}} 4E\left(\tau',s\right)$ and the second  to $ \frac{1}{2^s} 2E\left(2\tau',s\right)$. 
    
    \item Case $b_1$: We proceed in a similar way, and in this case we find that:
    
    \begin{equation*}
        \begin{split} S(\tau,s)&=\left(\sqrt{|G|}\right)^s\left(2\sum_{k,l}\frac{1}{\lambda(3k,l)^s}-\sum_{k,l}\frac{1}{\lambda(3k-1,l)^s}-\sum_{k,l}\frac{1}{\lambda(3k+1,l)^s}\right)\\ &=\left(\sqrt{|G|}\right)^s\left(3\sum_{k,l}\frac{1}{\lambda(3k,l)^s}-\sum_{k,l}\frac{1}{\lambda(k,l)^s}\right). 
        \end{split}
    \end{equation*}

    Rewriting in terms of $\tau'$ we get:
    
    $$S(\tau',s)=\left(\sqrt{|G|}\right)^s\left(3\frac{1}{3^s}E\left(\tau',s\right)-E\left(3\tau',s\right)\right)$$

\end{itemize}

All in all, carrying out the rest of the computations in the same fashion proves Theorem \ref{thm: global expression of combination of zeta functions}.

\subsection{The expression for the BCOV invariant}
\subsubsection{} We are now prepared to compute the BCOV invariant, renormalized as in the above sections, for all types of bielliptic surfaces. We recall that the various $\lambda$ are defined in Theorem \ref{thm: global expression of combination of zeta functions}. This expression, combined with \eqref{eq: quillen metric dist section bielliptic} directly proves Theorem \ref{Thm:F_1 bielliptic}.

\begin{theorem}\label{thm:computation of BCOV invariant}
    The BCOV invariant of bielliptic surfaces is given by the formula:
    $$\tau_{BCOV}=(2\im(\tau'))^2\prod_{j|N}\left(j^2|\Delta(j\tau')|^{\frac{1}{6}}\right)^{\lambda(j)}.$$
Moreover, $\prod_{j|N}j^{2\lambda(j)}=\frac{1}{N^4},\frac{1}{N^2},\frac{1}{N},1$ according to the cases $a,b,c,d$.

\end{theorem}
\begin{proof}
   From Lemma \ref{lemma:computation B(biell)} we already know $B(X)=\frac{1}{|G|}$. Moreover since $X$ carries a Ricci-flat metric and $\chi(X)=0$ we know $A(X)=1$. Hence, we only need to compute the spectral data.

    Recall from Theorem \ref{thm: global expression of combination of zeta functions} that we have the first equality of the following:
    \begin{equation}\label{eq:defpsi}
    \sum_{p,q}(-1)^{p+q}pq\zeta^{p,q}(s)=\frac{1}{K^s}\left(\sqrt{|G|}\right)^s\left(\sum_{j|N}\frac{1}{\left(\frac{N}{j}\right)^s}\lambda(j)E(j\tau',s)\right),
    \end{equation}
    where $G=\Z/N\Z\times \Z/K\Z$. The $\lambda(j)$ are explicit on each case and one has $\sum_{j|N}\lambda(j)=2$ and $\sum_{j|N}j\lambda(j)=0$. 

    Let us compute first the derivative of \eqref{eq:defpsi} at $s=0$ . One finds that it is given by: $$\left(-\log K + \log\left(\sqrt{|G|}\right)\right)\left(\sum_{j|N}\lambda(j)E(j\tau',0)\right)-\sum_{j|N}\left(\log\left(\frac{N}{j}\right)\lambda(j)E(j\tau',0)\right)+\sum_{k|N}\left(\lambda(j)E'(j\tau',0)\right)$$ 
    Using that $\sum\lambda(j)=2$ and, as noticed by the end of \ref{subsubsec:Kronecker} $E(j\tau',0)=-1$ for all $j$, combined with the Kronecker limit formula (see \ref{lemma:torsion ell.curve}), this expression simplifies to:
    \begin{equation}\label{derivative torsion combination split case}
        \begin{split}        
            2 \log K-\log|G|+2\log N-2\log(2)-2\log(\im(\tau'))-2\sum_{j|N}\lambda(j)\log(j)-\frac{1}{6}\sum_{j|N}\lambda(j)\log(|\Delta(j\tau')|)
        \end{split}
    \end{equation}

    Recalling that $\log K + \log N  = \log |G|$, the theorem follows by inserting \eqref{derivative torsion combination split case} into the definition  in \ref{eq:defbcov}. 
      
\end{proof}
\begin{remark}
    As we mentioned in the introduction it is easy to see that for surfaces $\mathcal{T}_{BCOV}=T(X,\Ocal_X)^{-2}$ hence the above is also a computation of analytic torsion for the trivial sheaf on bielliptic surfaces. 
\end{remark}

\section{$F_1^A$ of a bielliptic surface}\label{sect: F_1^A bielliptic}

The objective of this section is to compute the generating function of genus one Gromov--Witten invariants, defined in \eqref{Intro:def of F1A general}, for bielliptic surfaces. We follow the same conventions as in the previous sections and the implicit group in the definition is denoted by $G = \Z/K\Z \times \Z/N \Z$. The main result of this section is the following theorem, which is expressed in terms of the $\lambda(j)$ considered in Theorem \ref{Thm:F_1 bielliptic}.

\begin{theorem}\label{Thm:generating series GW computation}
    For all seven types of bielliptic surfaces,  we have: 
    $$F_1^A(\omega)=-\left(\sum_{j|N} \lambda(j)\log(\eta(q^j))\right)$$
\end{theorem}

Notice that this theorem together with Theorem \ref{Thm:F_1 bielliptic} immediately implies that the BCOV-mirror symmetry conjecture is true for bielliptic surfaces.

\subsection{Proof by degeneration}
The computation is done by means of the degeneration formula for Gromov--Witten invariants, precisely we apply \cite{LogDegenForm}. The main point of this formalism is to compute the invariants of a smooth member of a family of varieties in terms of logarithmic invariants of the singular fiber. Our computations are parallel to Blomme's computation \cite{blomme2025gromov} of all $F_g$ for $g\geq 2$.

\subsubsection{}From \cite[Subsect 2.4]{blomme2025gromov} we have a family $\mathcal{X}\to \D $ whose smooth fibers are bielliptic surfaces $X_t=\left(E\times F_t\right)/G$ for $t\neq 0$, and a single singular fiber $X_0$. This singular fiber is 
given as the union of two smooth irreducible components $X_1$ and $X_2$ glued along a divisor $D$. Both $X_1$ and $X_2$ are deformation equivalent to a product $\tilde{E}\times \mathbb{P}^1$ (so we treat them as this product), and $D$ has two connected components both isomorphic to $\tilde{E} = E/(\Z/K\Z)$. This singular fiber $X_0$ is an $\tilde{E}$-bundle over a Ner\'on $2$-gon. The monodromy action on a fiber $\tilde{E}$ is determined by a loop around the $2$-gon, which acts as multiplication by $\zeta_N=e^{\frac{2\pi i}{N} }$. Later on, in \ref{tab:monodromy} we recall Blomme's explicit computations of the monodromy action on $H^1$.

The aforementioned degeneration formula connects, by means of a count over certain graphs, the Gromov--Witten invariants $N_{1,\beta}$ of $X_t$ (which agree for all $t\neq 0$) with the relative/logarithmic invariants of $(X_1,D)$ and $(X_2,D)$, these count curves into the irreducible components of the singular fiber with prescribed tangency to the divisor.

\subsubsection{}
Before we start computing graph contributions let us remark that the degeneration formalism alone, gives us a lot of information. We begin by explaining a vanishing result for the relative invariants. We know (see \cite{blomme2025gromov} or \cite{SerranoBielliptics}) that the degree $2$ homology of the bielliptic  $X=(E\times F)/G$ is generated by two classes $\mathbf{e}$ and $\mathbf{f}$. Where $\mathbf{f}$ corresponds to $\frac{1}{K}$ times the class of the fiber of $X\to F/G$ and $\mathbf{e}$ corresponds to $\frac{1}{N}$ times the class of the fiber $X\to E/G$.
We have the following result, which is already mentioned in \cite[Subsect. 3.3.]{blomme2025gromov}. We include the result for completeness. 
\begin{proposition}\label{prop: vanishing Gw if not class ae}
The invariants $N_{1,a\mathbf{e}+b\mathbf{f}}$ vanish unless $b=0$
    
\end{proposition}
\begin{proof}
    The numerical version of the degeneration formula (see for example \cite[Theorem 1.6]{LogDegenForm}) 
    explicitly expresses the invariants $N_{1,\beta}$ as a sum of products of relative invariants.
    Hence, it is enough to show that the relative invariants vanish. 
    
    Notice that the variety $\tilde{E}\times \mathbb{P}^1$ inherits a $\C^\ast$-action from $\P^1$. Relative Gromov-Witten theory admits a torus localization formula.  Classically this is proven in \cite{GraberPhandaripandeLocalization} with two assumptions on the stack that are removed in \cite{AranhaetalVirtualLocrevisited}. The localization formula states that: 
    $$[\mathcal{M}_{g,n}((X,D),\beta)]^{virt}=\iota_\ast \left(\sum_i \frac{[\mathcal{M}_i^{\C^\ast}]^{virt}}{e(N_i^{virt})}\right)$$
    where $\mathcal{M}^{\C^\ast}_i$ runs over the connected components of the fixed point locus of the stack and $e(N_i^{virt})$ is the Euler class of its virtual normal bundle. We will finish the proof by showing that the fixed point locus is empty in our case. 
    
    Recall that we don't have point insertions for our curves, so in order for a curve $f:C\to \tilde{E}\times \mathbb{P}^1$ to lie on the fixed point locus, its image must be $\C^\ast$-equivariant. Since we have prescribed tangency conditions with $0$ and $\infty$ the curve cannot be entirely contained in $\tilde{E}\times 0$ or $\tilde{E}\times \infty$. 
    
    This means that a point of the form $(t,z)\in E\times \mathbb{P}^1$ with $z\in \C^\ast=\mathbb{P}^1\setminus\{0,\infty\}$ lies in the image of every curve. By properness and equivariance, the image of $C$ contains all of $t \times \P^1$. By connectedness of $C$ (see \cite{LogDegenForm} and dimension reasons it must be equal to it. Hence, if $b\neq 0$ the invariant locus of the moduli space is empty and so its virtual fundamental class vanishes.
    
\end{proof}

\subsubsection{} We now focus on computing the invariants for homology classes of the form $a\mathbf{e}$. In this context the result we apply is \cite[Theorem 1.5]{LogDegenForm} and we study curve counts with homology class $a\mathbf{e}$. We have, with notation introduced below:
\begin{equation}\label{eq:degenerationformula}
    N_{1,a\mathbf{e}}=\deg([\mathcal{M}_{1,0}(X,a\mathbf{e})]^{virt})=\sum_{\text{ certain graphs } \Gamma}\frac{\prod_{e \text{ edge }} \omega_e}{|\Aut \Gamma|}\deg(\Delta^!\prod_V[\mathcal{M}_V]^{virt}).
\end{equation}
Here $\Delta$ is the diagonal of a cartesian product of copies of the divisor indexed by the edges and $V$ denotes a vertex of the graph. We do not give the full definition of these graphs, see  \cite{LogDegenForm}, but mention that they are decorated graphs, with each vertex $V$ has an associated irreducible component $X_{r(V)}$ and genus $g_V$, and each edge $e$, by definition, has a tangency condition $\omega_e$.

 The graphs appearing in the formula decompose the topological type of the curves we are counting. That, is, they keep track of the fact that we are counting stable maps $f:C\to X_t$ with $C$ nodal of genus one and unpointed such that $f_\ast [C]=a\mathbf{e}$. This $\mathbf{e}$ splits inside $X_0$ as $\mathbf{e}=\mathbf{e}_1+\mathbf{e}_2$
where both $\mathbf{e}_i$ correspond to the fiber of $X_i=\tilde{E}\times \mathbb{P}^1\to \tilde{E}$

\subsubsection{}

Since the graphs parametrize curves into $X_0$ that intersect the divisor with the tangency condition it follows that these graphs must be cycles of $2h$ vertices, all with $g_v = 0$, corresponding to irreducible components of $C$. Moreover, all tangency conditions are equal, say to $w$, and we have that $a=h\cdot w$. The markings $r(V)\in \{1,2\}$ correspond to which component of $X_0$ the component of $C$ that the vertex represents is being mapped.  Consecutive vertices have alternating markings. 
 
Since $g_V=0$, the spaces $\mathcal{M}_V$ parametrize maps $f:\P^1\to X_{r(V)}$ with tangency condition $w$ at the divisor such that $f_\ast \P^1=w \mathbf{e}_{r(V)}$. The condition on homology means the projection onto the first factor $\widetilde{E}$ is constant. 
The projection to the second factor is a map $\P^1 \to \P^1$. The tangency conditions imply that we have a zero of order $w$ at 0 and a pole of order $w$ at $\infty$. Hence this map is equivalent to $z\mapsto z^w$. Its automorphism group is $\mu_w$. 

This means that 
\begin{displaymath}         \mathcal{M}_V=\tilde{E}_V\times B\mu_w
\end{displaymath} so $[\mathcal{M}_V]^{virt}=\frac{1}{w}[\tilde{E}_V]$. Since the graph is a loop the contribution of these $\frac{1}{w}$ cancels exactly the contribution of $\prod_e w$ in \eqref{eq:degenerationformula}.

The condition $a=h \cdot w$ determines each graph, so each divisor $h|a$ gives rise to such a graph. Each such graph has $h$ automorphisms, since the markings $r(V)$ must be respected. Hence, our formula \eqref{eq:degenerationformula} takes the form 
$$N_{1,a\mathbf{e}}=\sum_{h|a}\frac{\deg(\Delta^!\prod_V[\tilde{E}_V])}{h}=\sum_{h|a}\frac{A(h)}{h}.$$

The main result allowing explicit computation of $A(h)$ is the following:
\begin{theorem}\label{Thm:computation with traces}
We have $A(h)$ is the Lefschetz trace of $(M^h)^\ast$ on the cohomology of $\widetilde{E}$, in particular it is given by 
\begin{displaymath} 2-\operatorname{Tr}\left((M^h)^\ast|_{H^1}\right).
\end{displaymath}
Here $M$ is the monodromy given by one turn around the $2-gon$.
\end{theorem}
\begin{proof}
    We prove the result by intersection theoretic means. In order to compute the product of the diagonals we see that pull-back by the diagonal of the divisor $D$ can be reduced to a computation on $\tilde{E}$:
    
    Recall that $D=\tilde{E}_0\cup \tilde{E}_{\infty}\subset X_0$ where $E_0$ appears by gluing $D_1^-=\tilde{E}\times \{0\}$ on $\tilde{E}\times \P^1=X_1$ with $D^+_2=\tilde{E}\times \{\infty\}\subset\tilde{E}\times \P^1=X_2$ and $D_1^+=\tilde{E}\times\{\infty\}\subset X_1$ with $D_2^-=\tilde{E}\times \{0\}\subset X_2$.
    We do said identifications using multiplication by $\zeta_N$ on the $\tilde{E}$-factor.

    If we denote by $\varphi: D_1^-\to D_2^+$ and $\psi:D_2^-\to D_1^+$ the gluing maps which we interpret as maps from $\widetilde{E}$ to itself. We have, after identifying with $\widetilde{E}$, that $\psi\circ\varphi$ is the monodromy morphism $M$. 
    After identifications, the diagram where we take Gysin pull--back becomes:\\
    \begin{displaymath} \xymatrix{ & & \widetilde{E}^h \times \widetilde{E}^h \ar[d]^{\Delta_{\operatorname{Id}
     }\times \Delta_{\operatorname{Id}} }\\ 
    \widetilde{E}^{h}\times\widetilde{E}^h \ar[rr]^{G} & &\widetilde{E}^{2h} \times \widetilde{E}^{2h},}
    \end{displaymath}
where $G$ denotes the gluing map 
$(e,e')\mapsto (e,\phi(e),e',\psi(e'))$
for $e,e'$ $h$-tuples.

  This gives that  
    $$\Delta^!\left(\prod_V [E_V]\right)=G^!\left(\tilde{E}^h\times \tilde{E}^h \right)=\left(\Gamma_\phi \cdot\Gamma_\psi\right)^h.$$

Let us denote by $F:\tilde{E}\times\tilde{E}\to \tilde{E}\times \tilde{E}$ the map $(x,y)\mapsto (x,\phi(y))$. Then we have $\Gamma_\varphi =F^\ast\Delta_{\tilde{E}} $ and $F_\ast \Gamma_\psi=\Gamma_{M}$. Hence, since $F$ is an isomorphism, we obtain: 
\begin{equation}
    \begin{split}
      \deg(\Gamma_\phi \cdot \Gamma_\psi)=\deg(F_\ast(F^\ast\Delta_{\tilde{E}}\cdot \Gamma_\psi))=\deg(\Delta_{\tilde{E}}\cdot F_\ast \Gamma_\psi)=\deg(\Delta_{\tilde{E}}\cdot \Gamma_M)
    \end{split}
\end{equation}
from this we get $\deg(\Delta^!\prod_V[\tilde{E}_V])=\deg\left(\Delta_E\cdot \Gamma_{M^h}\right)$. We can now apply the Lefschetz fixed points formula (see e.g. \cite[Example 16.1.15]{FultonIntersection}) to obtain: 
\begin{equation}\label{eq: Coeffient of GW series}
    \begin{split}
    A(h)=\deg(\Delta^!\prod_V[E_V])=\deg(\Delta_E\cdot \Gamma_{M^h})=\sum_{j=0}^2(-1)^j\operatorname{Tr}\left((M^h)^\ast|_{H^j(E)} \right)=2-\operatorname{Tr}\left((M^h)^\ast|_{H^1}\right)
    \end{split}
\end{equation}
\end{proof}

\subsection{Computation of the monodromy factors. An example.}
Let us now compute $A(h)$ using Theorem \ref{Thm:computation with traces}. In order to do so, given the canonical basis of $H_1(\tilde{E},\Z)$, the monodromy action $M$ has the following shape (see \cite[Subsect 4.2.4]{blomme2025gromov}):
\begin{table}[h]\label{tab:monodromy}
\centering
\begin{tabular}{c c c c}
\textbf{a} & \textbf{b} & \textbf{c} & \textbf{d} \\
$\begin{bmatrix} -1 & 0 \\ 0 & -1 \end{bmatrix}$ &
$\begin{bmatrix} 0 & -1 \\ 1 & -1 \end{bmatrix}$ &
$\begin{bmatrix} 0 & -1 \\ 1 & 0 \end{bmatrix}$ &
$\begin{bmatrix} 0 & -1 \\ 1 & 1 \end{bmatrix}$ \\
\end{tabular}
\caption{different monodromy matrices for each type of surface}
\end{table}
\FloatBarrier

We know that $M^N=\id$ so we can write explicitly the different values of $A(h)$ for $h\in \Z/N\Z$ which are written in tuples in the following table:
\begin{equation}\label{eq:values-Ah}
\begin{array}{c c c c}
\mathbf{a} & \mathbf{b} & \mathbf{c} & \mathbf{d} \\
(0,4) &
(0,3,3) &
(0,2,4,2) &
(0,1,3,4,3,1)
\end{array}
\end{equation}
\subsubsection{}Let us now compute in full detail $F_1^A$ for type $a$. The other cases will be analogous. Based on \eqref{eq:values-Ah} we have that $F_1^A(q)=\sum_{a\geq 1}^\infty \sum_{h|a,\text{ odd }}\frac{4}{h}q^a$. With the convention that $\sigma_{-1}(x) = 0$ for $x \notin \N$ we have that: 
$$\sum_{h|a,\text{ odd }}\frac{1}{h}=\sum_{h|a}\frac{1}{h}-\sum_{h=2k, h|a}\frac{1}{h}=\sigma_{-1}(a)-\frac{1}{2}\sigma_{-1}\left(\frac{a}{2}\right).$$
 Putting everything together, and using that $-\log \eta(t) = -\frac{1}{24}\log  q+ \sum_a \sigma_{-1}(a) q^a$, we find that
$$F_1^A(q)=\sum_a\left(4\sigma_{-1}(a)-2\sigma_{-1}\left(\frac{a}{2}\right)\right)q^a=-4\log(\eta(\tau))+2\log(\eta(2\tau)).$$

    \newpage

\newpage
\printbibliography

\bigskip
\noindent
{\footnotesize
\textbf{Gabriel Abánades-Joglar}\\
Department of Mathematical Sciences,
Chalmers University of Technology and
University of Gothenburg, Sweden\\
\texttt{abanades@chalmers.se}

\medskip
\noindent
\textbf{Dennis Eriksson}\\
Department of Mathematical Sciences,
Chalmers University of Technology and
University of Gothenburg, Sweden\\
\texttt{dener@chalmers.se}
}
\end{document}